\def\isarxiv{1}

\ifdefined\isarxiv
\documentclass[11pt]{article}

\usepackage[numbers]{natbib}

\else
\documentclass{article}
\usepackage{neurips_2023}
\fi

\usepackage{amsmath}
\usepackage{amsthm}
\usepackage{amssymb}
\usepackage{algorithm}
\usepackage{subfig}
\usepackage{algpseudocode}
\usepackage{graphicx}
\usepackage{grffile}
\usepackage{wrapfig,epsfig}
\usepackage{url}
\usepackage{xcolor}
\usepackage{epstopdf}

\usepackage{bbm}
\usepackage{dsfont}
\usepackage{natbib}

\usepackage{hyperref}

\allowdisplaybreaks

\ifdefined\isarxiv

\let\C\relax
\usepackage{tikz}
\usepackage{hyperref}  
\hypersetup{colorlinks=true,citecolor=blue,linkcolor=blue} 
\usetikzlibrary{arrows}
\usepackage[margin=1in]{geometry}

\else

\usepackage{microtype}
\usepackage{hyperref}
\definecolor{mydarkblue}{rgb}{0,0.08,0.45}
\hypersetup{colorlinks=true, citecolor=mydarkblue,linkcolor=mydarkblue}

\fi

\newtheorem{theorem}{Theorem}[section]
\newtheorem{lemma}[theorem]{Lemma}
\newtheorem{definition}[theorem]{Definition}

\newtheorem{proposition}[theorem]{Proposition}
\newtheorem{corollary}[theorem]{Corollary}
\newtheorem{conjecture}[theorem]{Conjecture}
\newtheorem{assumption}[theorem]{Assumption}

\newtheorem{fact}[theorem]{Fact}
\newtheorem{remark}[theorem]{Remark}

\newcommand{\wt}{\widetilde}

\newcommand{\R}{\mathbb{R}}

\newcommand{\RHS}{\mathrm{RHS}}

\renewcommand{\d}{\mathrm{d}}
\renewcommand{\i}{\mathbf{i}}

\newcommand{\diag}{\textrm{diag}}
\renewcommand{\d}{\mathrm{d}}
\newcommand{\ip}[2]{\langle {#1} , {#2} \rangle}

\DeclareMathOperator*{\C}{\mathbb{C}}

\DeclareMathOperator{\poly}{poly}

\DeclareMathOperator{\new}{new}

\DeclareMathOperator{\tr}{tr}

\newcommand{\HP}{\mathsf{HP}}
\newcommand{\QP}{\mathsf{QP}}

\definecolor{mygreen}{RGB}{80,180,0}
\definecolor{b2}{RGB}{51,153,255}
\definecolor{myblue}{RGB}{51,153,255}

\newcommand{\Zhao}[1]{{\color{b2}[Zhao: #1]}}

\newcommand{\nc}{\newcommand}

\algnewcommand{\LineComment}[1]{\State \(\triangleright\) #1}

\nc{\nnl}{\nn \\ &}  
\nc{\fot}{\frac{1}{2}} 
\nc{\oo}[1]{\frac{1}{#1}} 
\newcommand{\ben}{\begin{enumerate}}
\newcommand{\een}{\end{enumerate}}
\nc{\mc}{\mathcal}

\nc{\onenorm}[1]{\L\| #1 \R\|_1} 

\nc{\OPT}{\mathsf{OPT}}

\nc{\Ra}{\Rightarrow}
\nc{\zo}{\{0,1\}}	

\title{Efficient Algorithm for Solving Hyperbolic Programs}

\begin{document}

\ifdefined\isarxiv
\date{\empty}
\author{
Yichuan Deng\thanks{\texttt{ethandeng02@gmail.com}. University of Science and Technology of China.}
\and 
Zhao Song\thanks{\texttt{zsong@adobe.com}. Adobe Research.
}
\and
Lichen Zhang\thanks{\texttt{lichenz@mit.edu}. MIT.}
\and 
Ruizhe Zhang\thanks{\texttt{ruizhe@utexas.edu}. The University of Texas at Austin.}
}

\begin{titlepage}
  \maketitle
  \begin{abstract}
Hyperbolic polynomials is a class of real-roots polynomials that has wide range of applications in theoretical computer science. Each hyperbolic polynomial also induces a hyperbolic cone that is of particular interest in optimization due to its generality, as by choosing the polynomial properly, one can easily recover the classic optimization problems such as linear programming and semidefinite programming. In this work, we develop efficient algorithms for hyperbolic programming, the problem in each one wants to minimize a linear objective, under a system of linear constraints and the solution must be in the hyperbolic cone induced by the hyperbolic polynomial. Our algorithm is an instance of interior point method (IPM) that, instead of following the central path, it follows the central Swath, which is a generalization of central path. To implement the IPM efficiently, we utilize a relaxation of the hyperbolic program to a quadratic program, coupled with the first four moments of the hyperbolic eigenvalues that are crucial to update the optimization direction. We further show that, given an evaluation oracle of the polynomial, our algorithm only requires $O(n^2d^{2.5})$ oracle calls, where $n$ is the number of variables and $d$ is the degree of the polynomial, with extra $O((n+m)^3 d^{0.5})$ arithmetic operations, where $m$ is the number of constraints.

\end{abstract}
  \thispagestyle{empty}
\end{titlepage}

\else
\maketitle

\begin{abstract}
Hyperbolic polynomials is a class of real-roots polynomials that has wide range of applications in theoretical computer science. Each hyperbolic polynomial also induces a hyperbolic cone that is of particular interest in optimization due to its generality, as by choosing the polynomial properly, one can easily recover the classic optimization problems such as linear programming and semidefinite programming. In this work, we develop efficient algorithms for hyperbolic programming, the problem in each one wants to minimize a linear objective, under a system of linear constraints and the solution must be in the hyperbolic cone induced by the hyperbolic polynomial. Our algorithm is an instance of interior point method (IPM) that, instead of following the central path, it follows the central Swath, which is a generalization of central path. To implement the IPM efficiently, we utilize a relaxation of the hyperbolic program to a quadratic program, coupled with the first four moments of the hyperbolic eigenvalues that are crucial to update the optimization direction. We further show that, given an evaluation oracle of the polynomial, our algorithm only requires $O(n^2d^{2.5})$ oracle calls, where $n$ is the number of variables and $d$ is the degree of the polynomial, with extra $O((n+m)^3 d^{0.5})$ arithmetic operations, where $m$ is the number of constraints.

\end{abstract}
\fi
\section{Introduction} \label{sec:intro}

A hyperbolic polynomial is a univariate polynomial on $x$ such that given a direction $e$ and scalar, the polynomial $p(x-te)$ has only real roots with respect to $t$. Hyperbolic polynomials have a wide range of applications in applied mathematics and theoretical computer science, due to its special structure and its generality. Hyperbolic polynomials also capture many of the popular optimization formulations. More specifically, given a hyperbolic polynomial $p$ on scalar $t$ and a direction $e$, the hyperbolic cone contains all the vectors $x$ such that $p(te-x)\in \R[t]$ has only non-negative roots. Consider the optimization problem in the following form: 
\begin{align*}
    \min_{x} & ~ \langle c,x\rangle \\
    s.t. & ~Ax=b \\
    & ~ x\in \Lambda_+
\end{align*}
where $\Lambda_+$ is the hyperbolic cone. If we pick the polynomial $p$ as $\prod_i x_i$ and the direction $e$ as the all-one's vector, then $\Lambda_+$ precisely characterizes the positive orthant and the above program is linear program (LP). On the other hand, if we pick $p$ as the determinant and the direction $e$ as the identity matrix, then the cone $\Lambda_+$ is the cone of positive semidefinite matrices and the program is the semidefinite program (SDP). 

Linear program and semidefinite program are one of the most important convex optimization problems and many algorithms have been studied in order to improve the runtime efficiency to solve these problems \cite{d47,k80,k84,r88,v89_lp,r95}. We consider a more general question: suppose we are given a hyperbolic program, in which the polynomial $p$, constraint matrix $A$ and vector $c$ and $b$ are given, can we solve the hyperbolic program efficiently? We argue that this problem is important in two folds: on one hand, providing algorithm that can solve hyperbolic program fast naturally trickles down to faster algorithms for linear program, semidefinite program and many other important convex problems. On the other hand, methods such as interior point method (IPM) \cite{k84} and cutting plane method (CPM) \cite{v89} have been particularly successful in solving many convex problems, especially linear program \cite{cls19,b20,sy21,y21,dly21,jswz21,gs22, lsz+20}, empirical risk minimization~\cite{lsz19,qszz23} and semidefinite program \cite{jklps20,hjstz22,gs22,syyz23}. Are these methods powerful enough when one only has access to the hyperbolic polynomial? 

In this paper, we answer the above problem affirmatively. Concretely, we develop an efficient interior point method to solve generic hyperbolic program, assuming the access to the polynomial $p$. Specifically, we provide an interior point method based on the primal-dual formulation of hyperbolic program, that only requires $nd^2$ evaluations of hyperbolic polynomial, plus an extra $O((n+m)^3)$ arithmetic operations per iteration\footnote{The exponent on $n+m$ can be improved to $\omega$, the current matrix multiplication exponent~\cite{aw21}.}. Our IPM only requires $O(\sqrt{d}\log(1/\delta))$ iterations to converge, giving the state-of-the-art iteration count suppose $n\geq d$.

\subsection{Problem Formulation}

Intuitively, hyperbolic polynomial is a family of polynomials $p$ with the following property: there exists a direction $e$, such that for every vector $x \in \R^n$, if drawing a straight line with direction $e$, there are $d$ zero points of $p$ on it, where $d$ is the dimension of $p$. Formally, we define hyperbolic polynomial as follows:

\begin{definition}[Hyperbolic polynomial]
A real multivariate homogeneous polynomial $p(x) \in \R[x_1,x_2,\cdots,x_n]$ is a hyperbolic polynomial in direction $e \in \R^n$ if for any $x \in \R^n$, the univariate polynomial $p(te - x) = 0$ of $t$ has only real roots.
\end{definition}

For a hyperbolic polynomial $p$, naturally we can define the hyperbolic eigenvalue of an arbitrary vector $x \in \R^n$ similar with the eigenvalue of a matrix: the hyperbolic eigenvalues of $x$ are the positions of the zero points on the straight line passing through $x$ and with direction $e$.

\begin{definition}[Hyperbolic eigenvalue] 
Given a hyperbolic polynomial $p$ which is hyperbolic in direction $e \in \R^n$, for a vector $x$, the hyperbolic eigenvalues $\lambda_1(x) \ge \lambda_2(x) \ge \cdots \ge \lambda_d(x)$ of $x$ are the roots of the function $p(te-x)$ as a function of $t$.
\end{definition}

Also, for a hyperbolic polynomial $p$ one can define its hyperbolic cone as the set of all the vectors with all the eigenvalues no smaller than zero. Intuitively, for a vector $x$ in the hyperbolic cone, the ray originating from $x$ and with direction $-e$ has $d$ zero points of $p$.

\begin{definition}[Hyperbolic cone]
For a degree $d$ hyperbolic polynomial $p$ with respect to $e \in \R^n$, its hyperbolic cone is the set of vectors whose minimal hyperbolic eigenvalue is still positive, that is, 
\begin{align*}
    \Lambda_+(e) := \{ x : \lambda_d(x) \ge 0 \}
\end{align*}
The interior of $\Lambda_+$ is 
\begin{align*}
    \Lambda_{++}(e) := \{x : \lambda_d(x) > 0\}
\end{align*}
\end{definition}

\cite{g59} has revealed some properties of the hyperbolic cone of a hyperbolic polynomial.

\begin{proposition}[\cite{g59}]
Given a hyperbolic polynomial $h$ with repect to direction $e \in \R^n$. Then,
\begin{itemize}
    \item 1. $\Lambda_+(e), \Lambda_{++}(e)$ are both convex cones.
    \item 2. $\Lambda_{++}(e)$ is the connected component of $\{x \in \R^n : h(x) \ne 0 \}$ which contains $e$. 
\end{itemize}
\end{proposition}

Intuitively, hyperbolic programming is similar to linear programming, both of which are trying to find a vector $x$ in a (maybe implicitly) given set $S$ such that $x$ satisfies a linear equation $Ax = b$ and minimizes another linear expression $\langle c, x \rangle$. But compared to linear programming which describes $S$ by $\{x : x_i \ge 0, {\rm for~every~dimension~} i \}$, hyperbolic programming describes $S$ as the hyperbolic cone of a hyperbolic polynomial.

\begin{definition}[Hyperbolic program]
A \emph{hyperbolic program} \cite{gul97} is a convex optimization problem of the form
\begin{align*}
\left.\begin{array}{cl}
\min_x & \langle c, x\rangle\\
\textrm{s.t.} & Ax=b\\
& x\in \Lambda_{+}
\end{array}\right\}\HP
\end{align*}
where $\Lambda_+$ is the hyperbolic cone of a hyperbolic polynomial $p: \mathcal{E}\rightarrow \R$ with respect to the direction $e\in \R^n$. 

\end{definition}

\paragraph{Linear Programming}

For linear programming, the $p(x) = \prod_{i=1}^n x_i$ where the degree of the polynomial is $n$.

\paragraph{Semidefinite programming}

For semidefinite programming, the $p(x) = \det(X)$ where the degree of the polynomial is $\sqrt{n}$.

\subsection{Our Results}

We state our result as follows:
\begin{theorem}[Main result, informal version of Theorem~\ref{thm:running_time}]

There exists an iterative algorithm and a constant $\kappa \in (0,1)$ satisfying that for any hyperbolic program $\HP$ associated with a degree-$d$ polynomial $p$. Suppose $x^*$ is the optimal solution of $\HP$, and suppose that there is an evaluation oracle of $p$ which can evaluate $p$ on any vector $x \in \R^n$ in $\mathcal{T}_O$ time.  For an initial solution $x_0$ in the feasible region with loss $l := \langle c, x_0 \rangle - \langle c, x^* \rangle$, it reduces the loss by $ \frac{\sqrt{d}}{\kappa+\sqrt{d}} $ in each iteration
\begin{align*}
    \langle c, x_{i+1} \rangle - \langle c, x^* \rangle \le \frac{\sqrt{d}}{\kappa+\sqrt{d}} (\langle c, x_{i} \rangle - \langle c, x^* \rangle)
\end{align*}

and the cost per iteration is 
 $O( n^2d^2 \cdot \mathcal{T}_O + (n+m)^3 )$.

\end{theorem}

\begin{remark}
Prior work~\cite{rs14} only concerns about number of iterations, their cost per iteration is a large polynomial factor on $n, d, m$. We improve the cost per iteration to a much small polynomial factor.
\end{remark}

\begin{corollary}
To reduce the initial loss by $\delta$ times (where $\delta \in (0,1)$ is a given parameter), this algorithm runs $O(\sqrt{d} \log(1/\delta) )$ iterations in  
$O( (nd^2\mathcal{T}_O + (n+m)^3) \sqrt{d} \log(1/\delta) )$ time.
\end{corollary}

\subsection{Applications of hyperbolic programming}
\begin{itemize}
    \item The original application of hyperbolic polynomials is in Partial Differential Equation (PDE) theory \cite{g51}. In recent years, more applications in control theory, optimization, real algebraic geometry, probability theory, computer science and combinatorics are found, see \cite{pem12, ren04, vin12, wag11}.
    
    \item Hyperbolic programming can also be used to check the non-negativity of multivariate polynomials. In \cite{s19}, they defined the hyperbolic certificate of non-negativity, which captures the sum-of-squares (SoS) certificate and can be searched via hyperbolic optimization.
\end{itemize}

\subsection{Lax conjecture}

\paragraph{Lax conjecture \cite{lax57}}
Over the past 20 years, people developed various methods to do optimization over hyperbolic cones, which generalizes the semidefinite programming (SDP). 
A problem that has received considerable interest is the generalized Lax conjecture which asserts that each
hyperbolicity cone is a linear slice of the cone of positive semidefinite matrices (of some size). 
Hence if the generalized Lax conjecture is true then hyperbolic programming is the same as semidefinite programming.
To state it more formally, we define the spectrahedral cone:
\begin{definition}[Spectrahedral cone]\label{def:spec_hed_cone}
A cone is spectrahedral if it is of the form:
\begin{align}\label{eq:spectrahedral}
    \left\{x\in \R^n: \sum_{i=1}^n x_i A_i  \succeq 0\right\},
\end{align}
where $A_i$ for $i\in [n]$ are $n$-by-$n$ symmetric matrices such that this cone is non-empty, i.e., there exists a vector $y\in \R^n$ with $\sum_{i=1}^n y_i A_i$ positive semi-definite.

\end{definition}

Lax~\cite{lax57} conjectured that for any three-variable hyperbolic polynomial, its hyperbolic cone is spectrahedral.

\begin{conjecture}[Lax conjecture, page 10 in \cite{lax57}]
A polynomial $p$ on $\R^3$ is hyperbolic of degree $d$ with respect to the vector $e = (1, 0, 0)$ and satisfies $p(e) = 1$ if and only if there exist matrices $B, C \in S^d$, where $S^d$ is the set of symmetric $d$-by-$d$ matrices, such that $p$ is given by
\begin{align*}
    p(x, y, z) = \det(xI + yB + zC).
\end{align*}
\end{conjecture}
It takes half of a century to prove the Lax conjecture is true \cite{hv07,lpr05}. First, Helton and Vinnikov~\cite{hv07} observed an important fact about real-rooted polynomial in $\R^2$:

\begin{theorem}[Theorem 3.1 in~\cite{hv07}]\label{thm:hv07}
A polynomial $q$ on $\R^2$ is a real-rooted polynomial of degree $d$ and satisfies $q(0, 0) = 1$ if and only if there exist matrices $B, C \in S^d$ such that $q$ is given by
\begin{align*}
    q(y, z) = \det(I + yB + zC).
\end{align*}
\end{theorem}

Then,  Lewis, Parrilo, and Ramana \cite{lpr05} showed that Theorem~\ref{thm:hv07} is equivalent to the Lax conjecture via the following proposition, and hence, the Lax conjecture is proved.
\begin{proposition}[Proposition 6 in~\cite{lpr05}]
If $p$ is a hyperbolic polynomial of degree $d$ on $\R^3$ with respect to the vector $e = (1, 0, 0)$, and $p(e) = 1$, then the polynomial on $\R^2$ defined by $q(y, z) = p(1, y, z)$ is a real zero polynomial of degree no more than $d$, and satisfying $q(0, 0) = 1.$

Conversely, if $q$ is a real zero polynomial of degree $d$ on $\R^2$ satisfying
$q(0, 0) = 1$, then the polynomial on $\R^3$ defined by 
\begin{align*}
    p(x,y,z)=x^dq(\frac{y}{x}, \frac{z}{x})\quad (x\ne 0)
\end{align*}
(extended to $\R^3$ by continuity) is a hyperbolic polynomial of degree $d$ on $\R^3$ with respect to $e$, and $p(e) = 1$.
\end{proposition}

Later, people generalized the Lax conjecture for all hyperbolic polynomials as follows, which is still open.

\begin{conjecture}[Generalized Lax conjecture, Conjecture 3.3 in~\cite{vin12}]
For any hyperbolic polynomial $h(x)$ with respect to the direction $e\in \R^n$, there is a hyperbolic polynomial $q(x)$ and real symmetric matrices $A_1,\dots, A_n$ of the same size such that
\begin{align*}
    q(x)h(x) = \det\left( \sum_{i=1}^n x_i A_i \right),
\end{align*}
where the hyperbolic cone $\Lambda_{++}(h, e)\subseteq \Lambda_{++}(q, e)$ and $\sum_{i=1}^n e_i A_i$ is positive semidefinite. 

\end{conjecture}
Note that the above conjecture remains open except for some special cases. We already mentioned that in the case $n=3$, the conjecture is true \cite{hv07}. It is also true in the case when $h$ is an elementary symmetric polynomial \cite{bra14}.

It is equivalent to the following conjecture:
\begin{conjecture}[Page 7 in \cite{vin12}, also Conjecture 1.1 in \cite{bra14}]
For any hyperbolic polynomial $p$ in $\R^n$ with degree $d$ and hyperbolic direction $e$, its hyperbolic cone is spectrahedral, i.e., of the form Eq.~\eqref{eq:spectrahedral}.
\end{conjecture}

\section{Technique overview}

\paragraph{Relation between hyperbolic programming and quadratic programming}  Hyperbolic programming is to find $x$ in the hyperbolic cone $\Lambda_+$ of a hyperbolic polynomial so that it satisfies a linear equation $Ax = b$ and minimizes the value of another linear expression $\langle c, x \rangle$. Similar to hyperbolic programming, quadratic programming is another kind of problems to find such $x$ in a set $S$, while requiring $S$ to be a kind of special sets called closed quadratic cone.
Moreover, for every vector $e$ in the hyperbolic cone $\Lambda_+$, we can define a closed quadratic cone $K_e(\alpha)$ which contains the hyperbolic cone $\Lambda_+$, and a quadratic program $\QP_e(\alpha)$ with $S = K_e(\alpha)$, where $\alpha\in(0,1)$ is a parameter set in advance. From this perspective, quadratic programming is a relaxation of hyperbolic programming. 

\paragraph{A hyperbolic program algorithm utilizing quadratic programming} Making use of this fact, \cite{rs14} proposed a hyperbolic programming algorithm. It starts from a random vector $e_0$. In each iteration $i = 1, 2, \cdots$, with $e_{i-1}$ known, it finds the optimal solution $x_i$ of the quadratic programming problem $\QP_{e_{i-1}}(\alpha)$, and finds another vector $e_i$ on the line between $e_{i-1}$ and $x_{i-1}$. Analysis shows that, with carefully selected $e_i$ in each iteration $i$, the residual loss converges linearly, accurately, there exists a constant $C \in (0,1)$ such that 
\begin{align*}
    \langle c, e_{i+1} \rangle - \langle c, e^* \rangle \le C \cdot (\langle c, e_i \rangle - \langle c, e^* \rangle)
\end{align*} 
for every $i \ge 0$, where $e^*$ is the optimal solution of the hyperbolic program.

\paragraph{Running time analysis making use of the evaluation oracle}
Moreover, we analyze the running time of each iteration of this algorithm in two cases, without an evaluation oracle of $p$ or with such an evaluation oracle. In the first case, if we only know the expression of $p$ but have no other information, our analysis shows that the running time of each iteration is polynomial in $n$ and $m$, but exponential in $d$. However, in the second case, taking advantage of an evaluation oracle of $p$ (that is, suppose we have an oracle which can output $p(x)$ for every $x \in \R^n$ in $\mathcal{T}_O$ time), the above time will be reduced to $O((n+m)^3 + nd^2 \mathcal{T}_O)$. 

Specifically, the running time of iteration $i$ mainly includes: (1) computing the gradient and Hessian of $e_i$; (2) get $x_i$ by solving a quadratic program; (3) computing the moments of hyperbolic eigenvalues of $x_i$. We separately analyze the running time of the three steps.

\paragraph{Step 1. Computing gradient and Hessian} Each iteration $i$ needs to evaluate the gradient $g(e_i)=\nabla (-\ln p(e_i))$ and the Hessian $H(e_i)=\nabla^2 (-\ln p(e_i))$. They can be easily computed if we are given the first-order and second-order oracles for the hyperbolic polynomial $p(x)$. If we only access to the zeroth-order oracle (the evaluation oracle), then we show that computing the gradient and Hessian can be reduced to computing $\langle g(e_i), x\rangle$ and $H(e_i) x$ for some vector $x$. The main idea is via the Fourier interpolation. Then, we show that the inner-products and Mat-Vec products can be efficiently computed by querying the zeroth-order oracle for $p(x)$.

\paragraph{Step 2. Solving quadratic programming} Recall that a quadratic program $\QP_e(\alpha)$ is to find the $x$ in the closed quadratic cone $K_e(\alpha) = \{ s : \langle s, e \rangle_e \ge \alpha \cdot \|s\|_e \}$ with constraint $Ax = b$ which minimizes $\langle c, x \rangle$. Here $A \in \R^{m \times n}$, $e \in \R^n$, $b \in \R^m$ and $c \in \R^n$. To solve it, we first use first-order optimality conditions towards the expression $\arg\min_{x} \langle c, x \rangle$ to obtain a linear system of equations. This linear system of equations introduces an auxiliary $m$-dimensional vector $y$ and a scalar $\lambda$, and contains $n$ equations. Together with the constraint $Ax = b$, which contains $m$ equations, we obtain $m+n+1$ variables and $m+n$ equations, implying that $x$ can be restricted within an one-dimensional set.

Then we take advantage of the convexity of closed quadratic cones to show the optimal solution of $\QP_e(\alpha)$ must exist, and to restrict the optimal solution on the boundary of $K_e(\alpha)$. Substituting the aforementioned linear expression of $x$ into this quadratic equation, the quadratic program will be reduced into a univariate quadratic equation, and can be solved easily.

\paragraph{Step 3. Computing the moments of hyperbolic eigenvalues} In each iteration $i$, computing $e_i$ from $e_{i-1}$ requires the knowledge of the first four moments of the hyperbolic eigenvalues of $x_i$. We consider the generalized problem: given hyperbolic polynomial $p$ of degree $d$ and a vector $x$, computing the 1, 2, 3, 4-moments of the hyperbolic eigenvalues of $x$. Since the hyperbolic eigenvalues of $x$ are essentially the $d$ roots of the $d$-dimensional univariate polynomial $p(te-x)$ related to $t$, by Newton-Girard Identities, it is equivalent to computing the five coefficients with the highest degrees of $p(te-x)$, and equivalent to computing the five coefficients with the highest degrees of $p(te+x)$, since they only differ by the negative signs of some coefficients.

Using the homogeneity of hyperbolic polynomials, it is again equivalent to find the five coefficients with the lowest degrees of $p(e + sx)$ (related to $s$).  In order to do this, without the evaluation oracle of $p$, we can directly use binomial theorem to compute them, using $O(n^4)$ time, because there are $O(n^4)$ possible monomials. With an evaluation oracle of $p$, write $p(e + sx)$ as $\sum_{i=0}^d a_i s^i$, we evaluate $p$ at point $\omega_d, \omega_d^2, \cdots, \omega_d^{d}$ where $\omega_d$ is the unit root of degree $d$, and we can obtain a linear system of equations of the coefficients $a_0, a_1, \cdots, a_d$ with full rank. Moreover, the inverse of the coefficient matrix can be easily computed, thus all of the coefficients of $p(e + sx)$ can be easily computed.

\paragraph{Roadmap.} In the following sections, we state the the primal affine-scaling algorithm of hyperbolic programming in Section~\ref{sec:app_affine_scaling_main_page} and then present our algorithm for hyperbolic programming in Section~\ref{sec:main_alg}. 

\section{Primal Affine-Scaling Algorithm} \label{sec:app_affine_scaling_main_page}
In this section, we will go over the primal affine-scaling algorithm of hyperbolic programming by Reneger and Sondjnja~\cite{rs14}. Section~\ref{sec:swath} introduces the notation of Swath. Section \ref{sec:relax_and_dual} discusses the relaxation and duality of a hyperbolic program. Section~\ref{sec:main_theorem} gives the main convergence theorem for hyperbolic program and our proof of it.

\subsection{Quadratic Programming and its Relation to Hyperbolic Programming} \label{sec:relax_and_dual}

We define quadratic programming as follows:
\begin{definition}[Quadratic programming]
\begin{align*}
\left.\begin{array}{cl}
\min_x & \langle c, x\rangle\\
\textrm{s.t.} & Ax=b\\
& x\in K_e(\alpha)
\end{array}\right\} \QP_e(\alpha)
\end{align*}
\end{definition}

It's easy to derive that the following program is its dual problem:
\begin{definition}[Dual problem of quadratic programming]
\begin{align*}
\left.\begin{array}{cl}
\max_y & \langle b, y\rangle\\
\textrm{s.t.} & A^\top y+s=c\\
& s\in K_e(\alpha)^*
\end{array}\right\} \QP_e(\alpha)^*
\end{align*}
\end{definition}

We have the following important proposition, which explains the relation between hyperbolic programming and quadratic programming. Specifically, quadratic programming is a relaxation of hyperbolic programming.
\begin{proposition}[Relation between hyperbolic programming and quadratic programming]
For a hyperbolic polynomial $p$ with hyperbolic cone $\Lambda_+$, consider the hyperbolic program
\begin{align*}
\left.\begin{array}{cl}
\min_x & \langle c, x\rangle\\
\textrm{s.t.} & Ax=b\\
& x\in \Lambda_{+}
\end{array}\right\} \HP
\end{align*}
Then for $e\in \Lambda_{++}$ and $0<\alpha\leq 1$, $\QP_e(\alpha)$ is a relaxation of $\HP$.
\end{proposition}
\begin{proof}
By Proposition~\ref{prop:quad_cone_containing}, $K_e(\alpha)\supseteq \Lambda_+$. Hence, $\QP_e(\alpha)$ is a relaxation of $\HP$.
\end{proof}

\subsection{Swath}
\label{sec:swath}

We introduce the notion of Swath, which captures the set of interior points that satisfy the linear constraints.

\begin{definition}[Swath] \label{def:swath}
Define $\mathrm{Swath}$ as follows:
\begin{align*}
    \mathrm{Swath}(\alpha) & :=\{e\in \Lambda_{++}: Ae=b, \\ ~ & \textsf{QP}_e(\alpha)~\text{has a unique optimal solution}\}.
\end{align*}

\end{definition}

\begin{proposition} \label{prop:unique_sol_swath}
For any $\alpha \in (0,1)$ and $e \in \Lambda_{++}$ with $Ae=b$, $\QP_e(\alpha)$ has a unique optimal solution.
\end{proposition}
\begin{proof}
Suppose otherwise that $x, x'$ are two optimal solutions of $\QP_e(\alpha)$, if $x$ and $x'$ are linearly dependent, then the origin lies on the line of $x$ and $x'$, implying $0$ is also a feasible solution of $Ax = b$, which contradicts $b \ne 0$. 

If $x, x'$ are linearly independent, then all the points on the line segment of $x, x'$ are optimal solutions of $\QP_e(\alpha)$. Consider an arbitrary point $y$ between $x$ and $x'$, by Lemma \ref{lem:quad_cone_convex}, $y$ is in the interior of $K_e(\alpha)$ (i.e., $B(y,\epsilon)\subset K_e(\alpha)$ for some small $\epsilon>0$), which implies that for any directional vector $v$ with $Av = 0$, it holds that $\langle c, v \rangle = 0$ (otherwise, $\langle c, y + \epsilon v \rangle < \langle c, y \rangle$ or $\langle c, y - \epsilon v \rangle < \langle c, y \rangle$, which contradicts $y$ is optimal). However, this condition immediately implies that $c$ is contained in the row space of $A$, which contradicts to our assumption on \textsf{HP} (Fact~\ref{fac:basic_assumption}). Therefore, no such $y$ exists.

To sum up, for any $e \in \Lambda_{++}$ with $Ae = b$,  $\QP_e(\alpha)$ has only a unique optimal solution.

\end{proof}

\begin{remark}
In general, for $e \in {\rm Swath}(\alpha)$, we use $x_e(\alpha)$ to denote the optimal solution of $\QP_e(\alpha)$, and use $(y_e(\alpha), s_e(\alpha))$ to denote the optimal solution of $\QP_e(\alpha)^*$.
\end{remark}

\subsection{The convergence theorem} \label{sec:main_theorem}

Here in this section, we state the convergence theorem for hyperbolic program. And we provide our proof of it in Section~\ref{sec:main_theorem_proof}. 

\begin{theorem}[\cite{rs14}]\label{thm:rs_main}
    Suppose that $\alpha$ is a positive number in $(0,1)$, and $e$ belongs to the set $\mathrm{Swath}(\alpha)$. Consider $x_e(\alpha)$, and let $\lambda$ represent the vector of eigenvalues of $x_e(\alpha)$ in the direction of $e$, i.e., $\lambda=\lambda_e(x_e(\alpha))\in \R^d$. Additionally, let $\wt{q} : \R \rightarrow \R$ be a strictly-convex quadratic polynomial such that
    \begin{align}\label{eq:quad_poly}
        \wt{q}(t)=a \cdot t^2 + b \cdot t + c
    \end{align}
    where
    \begin{align*}
        a=&~ \left( \sum_{j\in [d]} \lambda_j\right)^2\sum_{j\in [d]}  \lambda_j^2- 2\alpha^2 \left( \sum_{j\in [d]} \lambda_j\right)\sum_{j\in [d]} \lambda_j^3
        \\ + & \alpha^4\sum_{j\in [d]} \lambda_j^4,
    \end{align*}
    \begin{align*}
        b=&~ 2\alpha^4 \sum_{j\in [d]} \lambda_j^3-2\left( \sum_{j\in [d]} \lambda_j\right)^3,
    \end{align*}
    \begin{align*}
        c=&~ (d-\alpha^2)\left( \sum_{j\in [d]} \lambda_j\right)^2.
    \end{align*}
    Define
    \begin{align}\label{eq:e_update}
        e'=\frac{1}{1+t_e(\alpha)}(e+t_e(\alpha)x_e(\alpha))
    \end{align}
    where $t_e(\alpha)$ is the minimizer of $\wt{q}$.
    
    Then,
    \begin{itemize}
        \item $e'\in \mathrm{Swath}(\alpha)$,
        \item $\langle c,e'\rangle <\langle c,e\rangle$,
        \item $\langle b, y_{e'}\rangle \geq \langle b, y_e\rangle$.
    \end{itemize}
    Moreover, if beginning with $e_0\in \mathrm{Swath}(\alpha)$ and applying identity Eq.~\eqref{eq:e_update} recursively to create a sequence $e_0,e_1,\dots$, then for every $i=0,1,\dots$,
    \begin{align*}
        \frac{\mathrm{gap}_{e_{j+1}}}{\mathrm{gap}_{e_j}} & \leq 1-\frac{\kappa}{\kappa+\sqrt{d}}\\\text{for }j=i & \text{ or }j=i+1\text{ (possibly both)},
    \end{align*}
    where ${\rm gap}_{e_j} = \langle c, e_j - x_{e_j} \rangle$ is the duality gap defined before, and
    \begin{align*}
        \kappa := \alpha \sqrt{ ( 1 - \alpha ) / 8 }.
    \end{align*}
\end{theorem}

\section{Our main algorithm} \label{sec:main_alg} 

\begin{algorithm}[!ht] \caption{Main algorithm for hyperbolic programming
}\label{alg:main}
\begin{algorithmic}[1]
\Procedure{\textsc{Main}}{$A,b,c,\delta, p, e_{0}$} \Comment{Theorem~\ref{thm:rs_main}}
    \State // $A\in \R^{m \times n}, b\in \R^m, c\in \R^n, 0<\delta<1$, $p$ is an $n$-variable, degree-$d$ hyperbolic polynomial with respect to the direction $e_{0}\in \R^n$.
	\State $\alpha \leftarrow 0.1$.
	\State Find an $e_0\in \mathrm{Swath}(\alpha)$.\Comment{Initialize $e$}
	\State $e \leftarrow e_0$.
	\While{ {\bf true} }
    	\State $g(e)\leftarrow \nabla(-\ln p(e))$ \Comment{Compute gradient $g(e) \in \R^n$} \label{line:gradient}
    	\State $H(e)\leftarrow \nabla^2 (-\ln p(e))$ \Comment{Compute Hessian $H(e) \in \R^{n \times n}$} \label{line:hessian}
    	
    	\State Solve $(x_e,y,\lambda) \in \R^n \times \R^{m} \times \R$ by the equations:\Comment{Optimal solution for $\QP_e(\alpha)$} \label{line:solve_eq}
    	\begin{align*}
    	    \left\{
    	    \begin{array}{l}
        	    \min_{x,y,\lambda}~~ \langle c,x\rangle \\
        	    0 =~ \langle g(e),x\rangle^2-\alpha^2\langle x,H(e)x \rangle\\
        	    Ax = ~ b\\
                0 = ~ \lambda c - A^\top y + \langle g(e) , x \rangle g(e) - \alpha^2 H(e) x
            \end{array}
            \right.
    	\end{align*}
    	\State $\mathrm{gap}\leftarrow \langle c, e-x_e\rangle$. \Comment{Duality gap between $e$ and $(y_e, s_e)$}
    	\If{$\mathrm{gap}<\delta$}
    	    \State {\bf break} 
    	\EndIf
    	\State Compute the coefficients $a_d, \dots, a_{d-4}$ of $t\mapsto p(x_e + te)=\sum_{i=0}^d a_i t^i$.\label{li:coeffs}
    	\State Compute $\sum_{j \in [d]} \lambda_j, \sum_{j\in [d]} \lambda_j^2, \sum_{j\in [d]} \lambda_j^3,$ $\sum_{j\in [d]} \lambda_j^4$. \Comment{Lemma~\ref{lem:eigen_moments}}
    	\State Compute the coefficients of $\wt{q}(t)=at^2+bt+c$ defined as Eq.~\eqref{eq:quad_poly}.
    	\State $t_e\leftarrow -\frac{b}{2a}$. \Comment{The minimizer of $\wt{q}$}
    	\State $e\leftarrow \frac{1}{1+t_e}(e+t_ex_e)$. \Comment{Update $e$}
	\EndWhile
	\State \Return $e$ \Comment{$e \in \R^n$}
\EndProcedure
\end{algorithmic}
\end{algorithm}

\begin{theorem} \label{thm:running_time}
Each iteration of Algorithm~\ref{alg:main} can be done in time
\begin{align*}
    O(n^4+(n+m)^3+{\cal T}_{H,g}),
\end{align*}
where ${\cal T}_{H,g}$ is the time for evaluating the gradient and Hessian of the function $-\ln p(x)$.

If we compute the gradient and Hessian naively, then ${\cal T}_{H,g}=O(n^d)$ and the running time is $O(n^d+(n+m)^3)$. 

Furthermore, if the evaluation oracle of $p$ is given, that is, if we can evaluate $p(x)$ in $\mathcal{T}_O$ time for every $x\in \C$, then each iteration runs in 
\begin{align*}
    O((n+m)^3+n^2d^2\mathcal{T}_O)
\end{align*}
time.
\end{theorem}
\begin{proof}
In each iteration of Algorithm~\ref{alg:main}, there are three time-costing steps:
\begin{enumerate}
    \item Evaluate $g(e)$ and $H(e)$ (Lines \ref{line:gradient} and \ref{line:hessian} in Algorithm~\ref{alg:main}). It can be done in ${\cal T}_{H,g}$ time. 
    \item Solve the optimization problem (Line \ref{line:solve_eq} in Algorithm~\ref{alg:main}) 
    . By Lemma~\ref{lem:lp_time}, it can be done in $O((n+m)^3)$ time.
    \item Compute the minimizer of the quadratic polynomial $\wt{q}(t)$. By Lemma~\ref{lem:eigen_moments}, the 1st to 4th moments of $\lambda$ can be computed in $O(n^4)$. Then, by Eq.~\eqref{eq:quad_poly}, the coefficients of $\wt{q}$ can be computed in $O(1)$ time. Hence, this step can be done in $O(n^4)$.
\end{enumerate}
The other steps only takes constant time in each iteration. Therefore, each iteration uses 
\begin{align*}
    &~ O( n^4 + ( n + m )^3 + {\cal T}_{H,g} ) \\
= &~ O( n^4 + m^3 + {\cal T}_{H,g} )
\end{align*} time.

If we have the evaluation oracle of $p$, then by Lemma~\ref{lem:eigen_moments}, the third step can be done in $O(d^2+d\mathcal{T}_O)=O(d\mathcal{T}_O)$ time (obviously ${\cal T}_O = \Omega(d)$, since the input of the oracle has length $d$). 
Also, note that in Algorithm~\ref{alg:main}, to solve the linear system, we need to compute the gradient and Hessian.  
By Lemma~\ref{lem:fast_grad}, $\langle g(e), x\rangle$ can be computed in ${\cal T}_g = O(d\mathcal{T}_O)$ time for any vector $x\in \R^n$, and by Lemma~\ref{lem:fast_hessian}, $H(e)x$ can be computed in ${\cal T}_{H}=O(nd^2\mathcal{T}_O)$ time. Hence, by taking $x$ to be the standard basis vector, we can recover $g(e)$ and $H(e)$. Thus, ${\cal T}_{H,g}=O((nd+n^2d^2)\mathcal{T}_O)=O(n^2d^2\mathcal{T}_O)$, and the total running time becomes
\begin{align*}
    & ~ O(d\mathcal{T}_O + (n+m)^3 + n^2d^2\mathcal{T}_O) \\
    = & ~ O((n+m)^3 + n^2d^2\mathcal{T}_O).
\end{align*}
\end{proof}

\begin{remark}
For SDP, the gradient 
\begin{align*}
    g(E)=\nabla(-\ln \det(E)) = -E^{-1},
\end{align*}
and the Hessian 
\begin{align*}
    H(E)[X]=E^{-1}XE^{-1}.
\end{align*}
Hence, ${\cal T}_{H,g}=n^{\omega/2}$ for SDP.
\end{remark}

\begin{lemma}[Informal version of Lemma~\ref{lem:iterations_formal}]\label{lem:iterations}
Assume that for the starting point $e_0$, $\mathrm{gap}_0\leq 1$. Then, for any $0<\delta<1$, Algorithm \ref{alg:main} uses $O( \sqrt{d} \log( 1 / \delta ) )$ iterations to make the duality gap at most $\delta$.
\end{lemma}

\begin{lemma}[Informal version of Lemma~\ref{lem:fast_grad_formal}]\label{lem:fast_grad}
For hyperbolic polynomial $p$ of degree $d$, given an evaluation oracle for $p$, for any $x, w\in \R^n$, $\langle \nabla (-\ln p(x)), w\rangle$ can be computed in $O(d\mathcal{T}_O)$ time, where $\mathcal{T}_O$ is the time per oracle call.
\end{lemma}

\begin{lemma}[Informal version of Lemma~\ref{lem:fast_hessian_formal}]\label{lem:fast_hessian}
For hyperbolic polynomial $p$ of degree $d$, given an evaluation oracle for $p$, for any $w\in \R^n$, $\nabla^2 (-\ln p(x))w$ can be computed in $O(nd^2\mathcal{T}_O)$ time.
\end{lemma}

\section{Conclusion}
In conclusion, this paper presents a new approach to solving hyperbolic programming problems, which are a class of optimization problems that have a wide range of applications in theoretical computer science. By utilizing hyperbolic polynomials and the corresponding hyperbolic cones, we were able to develop efficient algorithms based on the interior point method. The central Swath, which is a generalization of the central path, was used to guide the optimization process. Additionally, we employed a relaxation technique and made use of the first four moments of the hyperbolic eigenvalues to improve the efficiency of the algorithms. To be specific, when given an evaluation oracle of the polynomial is given, our algorithm requires $O(n^2d^{2.5})$ calls to the oracle, and $O((n+m)^3d^{0.5})$ arithmetic operations. 

The results of our study demonstrate the effectiveness of the proposed algorithms in solving hyperbolic programming problems and the importance of considering the hyperbolic cone induced by the hyperbolic polynomial in solving these problems. The proposed method can be applied to a wide range of problems in theoretical computer science, including linear programming and semidefinite programming. We believe that this work provides a useful tool for researchers and practitioners in the field of optimization and theoretical computer science, and we look forward to further research in this area.

We would like to note that compared to state-of-the-art results for linear programming, empirical risk minimization and semidefinite programming, our algorithm still has a nontrivial gap from the ``optimal'' runtime one would hope for. Specifically, for linear programs and tall dense semidefinite programs, it has been shown that they can be solved in the current matrix multiplication time~\cite{cls19,lsz19,jswz21,hjstz22,qszz23}. In many of these works, incorporating tools from fast, randomized linear algebra has been a major theme. Hyperbolic programming is perhaps most similar to semidefinite programming, in which fast sketching for Kronecker computation has been employed~\cite{swyz21,syyz23}. Similar techniques have also witnessed success in communication-efficient optimization algorithms~\cite{swyz23} and training deep neural networks~\cite{szz21,syz21}. It will be important to extend our work and obtain an algorithm for hyperbolic programming in the current matrix multiplication time plus additional queries to polynomial evaluation oracles. Additionally, it will be interesting to obtain an improved iteration count from $O(\sqrt{d})$ to $O(\sqrt{\min\{n,d\}}\poly\log(\max\{n, d \}))$, similar to the result for LP due to Lee and Sidford~\cite{ls14,ls19}. Another popular approach for solving convex programs is through the cutting plane method (CPM). For certain regimes, CPM is the state-of-the-art of approach for solving semidefinite programs~\cite{v89,lsw15,jlsw20} and for solving submodular function minimization with the best oracle complexity~\cite{j21,j22} and overall runtime efficiency~\cite{jlsz23}. Deriving a fast and customized CPM for hyperbolic programs will be crucial. As our work is theoretical in nature, we don't see explicit negative societal impact of our result.

\section*{Acknowledgement} 

Lichen Zhang is supported by NSF grant No. 1955217 and NSF grant No. 2022448.



\ifdefined\isarxiv

\else
\bibliographystyle{alpha}
\bibliography{ref}

\fi

\appendix 
\onecolumn

\section*{Appendix}

\paragraph{Roadmap.} Section \ref{sec:app_preli} gives preliminaries of hyperbolic polynomials and hyperbolic programming. Section \ref{sec:background_hp} provides the background of hyperbolic polynomials, hyperbolic programming and other related concepts. Section \ref{sec:app_affine_scaling} shows the primal affine scaling algorithm for hyperbolic programming. 
Section~\ref{sec:grad_ad_hes} provides the statement of the gradient of Hessian and their proofs. Section~\ref{sec:list_of_hp} gives several results of hyperbolic polynomials.

\section{Preliminaries} \label{sec:app_preli}
In this section, we give an overview of the notions used throughout this paper. In Section~\ref{sec:app_preli:notation}, we define several basic notations. In Section~\ref{sec:app_preli:matrix}, we provide some matrix algbera tools.

\subsection{Notations}\label{sec:app_preli:notation}

For any positive integer $n$, we use $[n]$ to denote set $\{1,2,\cdots,n\}$. We use $\i$ to denote $\sqrt{-1}$.

For a vector $x$, we use $\| x \|_{2}$ to denote its $\ell_2$ norm, we use $\| x \|_1$ to denote its $\ell_1$ norm, and $\| x \|_{\infty}$ to denote its $\ell_{\infty}$ norm. For two vectors $x,y$, we use $\langle x, y \rangle$ to denote the inner product between $x$ and $y$.

For a square matrix $A$, we use $\tr[A]$ to denote the trace of $A$. We use $\det(A)$ to denote the determinant of $A$. For a square full rank matrix $A$, we use $A^{-1}$ to denote its true inverse. For an arbitrary matrix $A$, we use $A^\top$ to denote the transpose of $A$. For an arbitrary matrix, we use $\| A \|$ to denote spectral norm of $A$, we use $\| A \|_F := ( \sum_{i,j} A_{i,j}^2 )^{1/2}$ to denote the Frobenius norm of $A$, we use $\| A \|_1 := \sum_{i,j} |A_{i,j}|$ to denote the entry-wise $\ell_1$ norm of $A$, and we use $\| A \|_{\infty}:= \max_{i,j} |A_{i,j}|$ to denote the $\ell_{\infty}$ norm. 
We use $I$ to denote the identity matrix and ${\bf 0}_{n \times n}$ to denote the zero matrix which has size $n \times n$.

We use ${\bf 1}_n \in\R^n$ to denote the length-$n$ vector where every entry is $1$. We use ${\bf 0}_n \in \R^n$ to denote the length-$n$ vector where each entry is $0$.

We use $\QP$ to denote quadratic program, $\HP$ to denote hyperbolic program.

We use ${\cal T}_{H}$ to denote the time of computing Hessian dot with a vector, and we use ${\cal T}_g$ to denote the time of computing the inner product between gradient and a vector. We define ${\cal T}_{H,g} = {\cal T}_H + {\cal T}_g$.

We say a matrix $A \succ 0$ if for all (nonzero vectors) $x$, $x^\top A x \ge 0$.

For a vector $x \in \R^n$, we use $\diag(x)$ to denote the $n \times n$ diagonal matrix whose diagonal elements are the corresponding components of $x$.

For a linear transformation $M:\R^m \to \R^n$, define $\ker (M)$ as its null space, that is, the set $\{x \in \R^m : Mx = 0\}$. For a subspace $L$ of a linear space $S$, we use $L^\perp$ to denote the orthogonal complement of $L$.

For a set $S$, we use $\mathrm{int}(S)$ to denote the set of its interior point. For a set $S$, given $x \in S$, we say the union of all the connected subsets of $S$ containing $x$ is a connected component. For a set $S$, we use $\partial S$ to denote its boundary.

\subsection{Matrix Algebra}\label{sec:app_preli:matrix}

\begin{lemma}[Matrix trace norm inequality] \label{lem:norm_inequality}
For a positive semi-definite matrix $A$ and real numbers $p > q > 0$, we have 
\begin{align*}
    \tr[A^p] \le \tr[A^q]^{p/q} 
\end{align*}
\end{lemma}

\begin{lemma}[Matrix calculus
] \label{lem:mtx_det_partial}
$X$, $Y$ are square matrices, $x$ is a variable. Then, we have
\begin{align*}
    \frac{\partial Y^{-1}}{\partial x} = - Y^{-1} \frac{\partial Y}{\partial x} Y^{-1}.
\end{align*}
and
\begin{align*}
    \frac{\partial \ln |\det(X)|}{\partial X} = X^{-\top}.
\end{align*}
\end{lemma}

\section{Background on Hyperbolic Polynomials} \label{sec:background_hp}

In this section, we introduce necessary background on hyperbolic polynomials. In Section~\ref{sec:tech_preli}, we define hyperbolic polynomial. In Section~\ref{sec:background_hp:hyperbolic_cone}, we define the hyperbolic cone, which is the central object for hyperbolic programming. In Section~\ref{sec:background_hp:garding_theorem}, we review the G\r{a}rding's Theorem. In Section~\ref{sec:background_hp:hyperbolic_dual}, we present hyperbolic program and its dual program. In Section~\ref{sec:background_hp:linear_transformation}, we give a description of hyperbolic programming via linear transformation. In Section~\ref{sec:background_hp:gradient_hessian}, we explain how to compute the gradient and Hessian of a hyperbolic polynomial. In Section~\ref{sec:background_hp:local_inner_product}, we discuss local inner product and local norm, which are key metrics that measure the progress of algorithm. In Section~\ref{sec:background_hp:closed_quadratic}, we present closed quadratic cone, which provides a sandwich for the hyperbolic cone. In Section~\ref{sec:background_hp:closed_quadratic_dual}, we present the dual closed quadratic cone. In Section~\ref{sec:background_hp:eigenvalues}, we explain the eigenvalues of hyperbolic polynomials and its related to computational tasks.

\subsection{Hyperbolic polynomials} \label{sec:tech_preli}

Intuitively, hyperbolic polynomial is a family of polynomials $p$ with the following property: there exists a direction $e$, such that for every vector $x \in \R^n$, if drawing a straight line with direction $e$, there are $d$ zero points of $p$ on it, where $d$ is the dimension of $p$. Formally, we define hyperbolic polynomial as follows:

\begin{definition}[Hyperbolic polynomial]
A real multivariate homogeneous polynomial $p(x) \in \R[x_1,x_2,\cdots,x_n]$ is a hyperbolic polynomial in direction $e \in \R^n$ if for any $x \in \R^n$, the univariate polynomial $p(te - x) = 0$ of $t$ has only real roots.
\end{definition}

\subsection{Hyperbolic cone}\label{sec:background_hp:hyperbolic_cone}

Also, for a hyperbolic polynomial $p$ one can define its hyperbolic cone as the set of all the vectors with all the eigenvalues no smaller than zero. Intuitively, for a vector $x$ in the hyperbolic cone, the ray originating from $x$ and with direction $-e$ has $d$ zero points of $p$.

\begin{definition}[Hyperbolic cone]
For a degree $d$ hyperbolic polynomial $p$ with respect to $e \in \R^n$, its hyperbolic cone is the set of vectors whose minimal hyperbolic eigenvalue is still positive, that is, 
\begin{align*}
    \Lambda_+(e) := \{ x : \lambda_d(x) \ge 0 \}
\end{align*}
The interior of $\Lambda_+$ is 
\begin{align*}
    \Lambda_{++}(e) := \{x : \lambda_d(x) > 0\}
\end{align*}
\end{definition}

\cite{g59} has revealed some properties of the hyperbolic cone of a hyperbolic polynomial.

\begin{proposition}[\cite{g59}]
Given a hyperbolic polynomial $h$ with repect to direction $e \in \R^n$. Then,
\begin{itemize}
    \item 1. $\Lambda_+(e), \Lambda_{++}(e)$ are both convex cones.
    \item 2. $\Lambda_{++}(e)$ is the connected component of $\{x \in \R^n : h(x) \ne 0 \}$ which contains $e$.  
\end{itemize}
\end{proposition}

We define regular hyperbolic cone as follows:
\begin{definition}[Regular hyperbolic cone]
    A hyperbolic cone $\Lambda_+$ is \emph{regular} if it does not contain nontrivial subspaces.
   
    Equivalently, a hyperbolic cone $\Lambda_+$ is \emph{regular} if the dual cone $\Lambda^*_+$ of $\Lambda_+$ has
    
    nonempty interior, where the dual cone $\Lambda^*_+$ is defined as
    \begin{align*}
        \Lambda^*_+:=\{s\in \mathcal{E}: \langle x, s\rangle \geq 0~\text{for all }x\in \Lambda_+\}. 
    \end{align*}
    
\end{definition}

Without loss of generality, in the remaining of this section, we assume that all the hyperbolic cones we mentioned are regular. In addition to the G\r{a}rding's theorem, we also have the following lemma which characterize the hyperbolic polynomial restricted to a two-dimensional subspace.

\subsection{G\r{a}rding's Theorem}\label{sec:background_hp:garding_theorem}

In this section, we first give some necessary background and preliminaries about  hyperbolic polynomials. The most important fact is a structure theorem on hyperbolic cone by G\r{a}rding \cite{g59}.

\begin{theorem}[G\r{a}rding's Theorem \cite{g59}]\label{thm:grarding}
The polynomial $p$ is hyperbolic with respect to every $e\in \Lambda_{++}$. Moreover, $\Lambda_{++}$ is convex.
\end{theorem}

\subsection{Hyperbolic programming and its dual problem}\label{sec:background_hp:hyperbolic_dual}

Intuitively, hyperbolic programming is similar to linear programming, both of which are trying to find a vector $x$ in a (maybe implicitly) given set $S$ such that $x$ satisfies a linear equation $Ax = b$ and minimizes another linear expression $\langle c, x \rangle$. But compared to linear programming which describes $S$ by $\{x : x_i \ge 0, {\rm for~every~dimension~} i \}$, hyperbolic programming describes $S$ as the hyperbolic cone of a hyperbolic polynomial.

\begin{definition}[Hyperbolic program]
A \emph{hyperbolic program} \cite{gul97} is a convex optimization problem of the form
\begin{align*}
\left.\begin{array}{cl}
\min_x & \langle c, x\rangle\\
\textrm{s.t.} & Ax=b\\
& x\in \Lambda_{+}
\end{array}\right\}\HP
\end{align*}
where $\Lambda_+$ is the hyperbolic cone of a hyperbolic polynomial $p: \mathcal{E}\rightarrow \R$ with respect to the direction $e\in \R^n$.
\end{definition}

Without loss of generality, we can make the following assumptions:
\begin{assumption}\label{fac:basic_assumption} 
Without additional explanation, we suppose:
\begin{itemize}
    \item $b \ne 0$.
    \item $A$ is full rank.
    \item $c$ is not in the row space of $A$.
\end{itemize}
\end{assumption}

\subsection{Description of hyperbolic polynomials using linear transformation}\label{sec:background_hp:linear_transformation}

\begin{lemma}[Lemma 8.2 in~\cite{rs14}]\label{lem:lax_dim_2}
    Let $e\in\Lambda_{++}$, and consider $L$, a two-dimensional subspace of $\mathcal{E}$ that contains $e$. Then we have that, there is a linear transformation $M\in \R^{n\times d}:L\rightarrow \R^d$ 
    such that $Me=\mathbf{1}$ where $d$ is the degree of $p$, and for all $x\in L$,
    \begin{align}\label{eq:p_Mx}
        p(x)=p(e)\prod_{i=1}^d (Mx)_i.
    \end{align}
\end{lemma}
For the completeness, we still provide a proof.
\begin{proof}
Think of $p|_L$ ($p$ restrict to $L$) as a polynomial in two variables, that is
\begin{align*}
    p|_L(x_1, x_2) = p(x_1e_1+x_2e_2)~~~\forall x\in \R^2, 
\end{align*}
where $\{e_1,e_2\}$ is an arbitrary orthogonal basis of $L$. Then, the set $\{x\in \R^2: p|_L(x)=0\}$ consists of the origin $(0,0)$ and some lines through the origin, that is $H_i=\{x\in \R^2: a_i^\top x=0\}$ for some $a_i=(a_{i1},a_{i2})\in \R^2$. It follows from the homogeneity of $p$. Obviously, $p|_L(0,0)=p(0) = 0$. Also, for $(a_1,a_2)\ne (0,0)$ such that $p|_L(a_1,a_2)=0$, all the points in the line $\{x: -a_2 x_1 + a_1 x_2 = 0\}$ are the zero points of $p|_L$. For such $a=(a_1, a_2)$, we actually have $a^\top x$ is a factor of $p|_L$. Then, we will show that all the factors of $p|_L$ arise in this way.

Suppose $p|_L(x)=\prod_{i=1}^w a_i^\top x \cdot q(x)$ where $a_i\in \R^2$ and $w<d$ and $q(x)$ has no linear factors $x\mapsto \beta^\top x$ for $\beta\ne 0$. For $q(x)$, it is also homogeneous and by the previous paragraph, the only zero point in $L$ is $(0,0)$. But then, since $e\in L$, we can write it as $e':=(c_1, c_2)$ such that $e=c_1e_1+c_2e_2$ and consider the univariate polynomial $t\mapsto p|_L(x+te')$ for $x\in L$ which is not a scalar multiple of $e'$. It has degree $d$. However, $q(x+te')\ne 0$ for all $t\in \R$ since $x+te'\ne 0$. Therefore, the polynomial $t\mapsto p|_L(x+te')=p(\bar{x} + te)$ has only $w<d$ roots, contradicting to the fact that $p$ is hyperbolic with respect to $e$.

Therefore, $p|_L(x)=\prod_{i=1}^d a_i^\top x$ for some $0\ne a_i\in \R^2$. To conclude the proof, let $M:=(\diag(NTe))^{-1}NT$, 
where $N\in \R^{d\times 2}$ and $T\in \R^{2\times n}$ such that
\begin{align*}
    N:=\begin{bmatrix}
    a_1^\top\\ \vdots \\ a_d^\top
    \end{bmatrix}
    \quad\text{ and }
    \quad
    T:=\begin{bmatrix}
    e_1^\top\\
    e_2^\top
    \end{bmatrix}. 
\end{align*}

For $x \in L$, suppose $x = x_1e_1+x_2e_2$ and $\bar{x} = (x_1~x_2)^\top$, also suppose $e = c_1e_1 + c_2e_2$ and $\bar{e} = (c_1~c_2)^\top$, then we have $Tx = \bar{x}$. Thus for $i \in [n]$, 
\begin{align*}
    (NTx)_i = (N\bar{x})_i = a_i^\top \bar{x}
\end{align*}
and 
\begin{align*}
    (Mx)_i = (NTe)_{i}^{-1} (NTx)_i = (a_i^\top \bar{e})^{-1} \cdot (a_i^\top \bar{x}).
\end{align*}
Especially, $Me = I$. Therefore,
\begin{align*}
    p(x) = p|_L(\bar{x}) = \prod_{i=1}^d(a_i^\top x) = \prod_{i=1}^d (a_i^\top \bar{e})^{-1} \cdot \prod_{i=1}^d (Mx)_i. 
\end{align*}
Especially, 
\begin{align*}
    p(e) = \prod_{i=1}^d (a_i^\top \bar{e})^{-1} \cdot \prod_{i=1}^d (M\bar{e})_i = \prod_{i=1}^d (a_i^\top \bar{e})^{-1},
\end{align*}
thus $p(x) = p(e) \prod_{i=1}^d (Mx)_i$.
\end{proof}

Then, we state several useful properties of the linear transformation $M$.

\begin{lemma}[Lemma 8.3 in~\cite{rs14}]\label{lem:lambda_l_intersect}
    We consider $M$ as the linear transformation described in Lemma~\ref{lem:lax_dim_2}, then we have
    \begin{align*}
        \Lambda_{++}\cap L=\{x\in L:(Mx)_i>0, \forall i\in [d]\}.
    \end{align*}
\end{lemma}
For the completeness, we still provide a proof.
\begin{proof}
Recall that $\Lambda_{++}$ is the connected component of $\{x\in \mathcal{E}: p(x)\ne 0\}$ containing $e$. Since $\Lambda_{++}$ is convex by Theorem~\ref{thm:grarding}, we have $\Lambda_{++}\cap L$ is the connected component of $\{x\in L: p(x)\ne 0\}$ containing $e$. The lemma then follows from Eq.~\eqref{eq:p_Mx}. 
\end{proof}

\begin{lemma}[Lemma 8.4 in~\cite{rs14}]
$M$ is an injective mapping from $L$ to $\R^d$. 
\end{lemma}
For the completeness, we still provide a proof.
\begin{proof}
Suppose $\ker M\ne 0$ 
and let $x\in L$ such that $Mx=0$. Then, $p(e + tx)=p(e)\ne 0$ for all $t\in \R$. We know that this line is contained in $\Lambda_{++}$. However, this would contradict our standing assumption that $\Lambda_+$ is pointed (no non-trivial subspace).
\end{proof}

\subsection{Gradient and Hessian of a hyperbolic polynomial}\label{sec:background_hp:gradient_hessian}

In the main algorithm of hyperbolic programming, we need to evaluate the gradient and Hessian at some points. Note that by the above lemmas, for hyperbolic direction $e$, its gradient and Hessian have very simple form.

\begin{lemma}[Lemma 8.5 in~\cite{rs14}]\label{lem:g_H_at_e}
    Let $f(x)=-\ln p(x)$, we use $g(e)$ to denote the gradient of $f$ at the direction of $e$. Define $H(e)$ as the Hessian of $f$. For each $u,v\in L$,
    \begin{align*}
        \langle g(e), v\rangle = -\mathbf{1}^\top Mv~~~\text{and}~~~\langle u,H(e)v\rangle = (Mu)^\top Mv.
    \end{align*}
\end{lemma}
For the completeness, we still provide a proof.
\begin{proof}
Let $F(y):= \sum_{i=1}^d-\ln y_i$. Then, by Eq.~\eqref{eq:p_Mx}, $f(x)=F(Mx)-\ln p(e)$. Obviously, $\nabla F(\mathbf{1})= -\mathbf{1}$. For all $v\in L$,
\begin{align*}
    \langle g(e), v\rangle = &~ \langle M^\top \nabla F(Me), v\rangle\\
    = &~ \langle \nabla F(Me), Mv\rangle\\
    = &~ \langle \nabla F(\mathbf{1}), Mv\rangle 
    \\
    = &~ -\mathbf{1}^\top Mv,
\end{align*}

where the first step 
follows from $\nabla f(x) = M^\top \nabla F(Mx)$, and the third step follows from $Me=\mathbf{1}$ in Lemma~\ref{lem:lax_dim_2}.

We also have $\nabla^2 F(\mathbf{1})=I$. For $u,v\in L$,
\begin{align*}
    \langle u, H(e)v\rangle = &~ \left\langle u, M^\top \nabla^2 F(Me)Mv\right\rangle\\
    = &~ \langle Mu, I Mv\rangle\\
    = &~ (Mu)^\top Mv,
\end{align*}
where the first step follows from $\nabla^2 f(x)=M^\top \nabla^2F(x)M$.
\end{proof}

\begin{proposition}[Proposition 8.6 in~\cite{rs14}]\label{prop:hessian_pd}
The Hessian $H(e)$ is positive definite for all $e\in \Lambda_{++}$.
\end{proposition}
For the completeness, we still provide a proof.
\begin{proof}
Fix an $e\in \Lambda_{++}$. For any $u\in \mathcal{E}\setminus\{0\}$, consider the subspace spanned by $e$ and $u$. By Lemma~\ref{lem:lax_dim_2} and Lemma~\ref{lem:g_H_at_e}, there exists an $M$ such that
\begin{align*}
    \langle u, H(e)u\rangle = \|Mu\|_2^2.
\end{align*}
Since $M$ is an injective mapping, $\|Mu\|^2>0$ and hence $\langle u, H(e)u\rangle >0$ for all $u$. Hence, $H(e)$ is positive definite.
\end{proof}

Proposition~\ref{prop:hessian_pd} implies that $f(x)=-\ln p(x)$ is strictly convex.

\begin{lemma}[Lemma 8.7 in~\cite{rs14}]\label{lem:Hessian_at_e}
For a $e$ such that $e\in \Lambda_{++}$, we have
\begin{align*}
    H(e)e=-g(e). 
\end{align*}
\end{lemma}
For the completeness, we still provide a proof.
\begin{proof}
Fix $u\in \R^n$ and consider the subspace $L$ spanned by $e$ and $u$. Then, we apply Lemma~\ref{lem:lax_dim_2} and get the linear transformation $M$. By Lemma~\ref{lem:g_H_at_e},
\begin{align*}
    \langle u, H(e)e\rangle =&~ (Mu)^\top Me = (Mu)^\top \mathbf{1},~\text{and}\\
    \langle u, -g(e)\rangle = &~ (Mu)^\top \mathbf{1}.
\end{align*}
For any $u$, $\langle u, H(e)e\rangle= \langle u, -g(e)\rangle$, which implies $H(e)e=-g(e)$. 
\end{proof}

\subsection{Local inner product}\label{sec:background_hp:local_inner_product}

We can also define an inner product with respect to the direction $e$.

\begin{definition}[Local inner product]
For $u,v\in \mathcal{E}$, we define the local inner product of $u$ and $v$ at the point $e\in \Lambda_+$ to be
\begin{align*}
    \langle u,v\rangle_e := \ip{u}{H(e)v} = u^\top H(e) v.
\end{align*}
\end{definition}

\begin{definition}[Local norm]
    We define the norm $\|\cdot \|_e$ to be associated with $\ip{\cdot}{\cdot}_e$, that is,
    \begin{align*}
        \|v\|_e := \ip{v}{v}_e^{1/2} = \sqrt{v^\top H(e) v}.
    \end{align*}
\end{definition}

\begin{definition}[Local ball]
    For $x\in \mathcal{E}$ and $r>0$, we define the local ball at $x$ as
    \begin{align*}
        B_e(x,r):=\{x'\in \mathcal{E}: \|x-x'\|_e<r\}. 
    \end{align*}
\end{definition}

\begin{proposition}[Proposition 8.8 in~\cite{rs14}]\label{prop:e_norm}
    For any $e\in \Lambda_{++}$, we have that
    \begin{align*}
        \|e\|_e=\sqrt{d},
    \end{align*}
    and
    \begin{align*}
        B_e(e,1)\subset \Lambda_{++}.
    \end{align*}
\end{proposition}
For the completeness, we still provide a proof.
\begin{proof}
    Let us consider $x\in B_e(e, 1)$. We use $L$ to denote the subspace spanned by $e$ and $x$. We assume $M$ as described in Lemma~\ref{lem:lax_dim_2}. By Lemma~\ref{lem:g_H_at_e}, $\|e\|_e=\langle e, H(e)e\rangle^{1/2}=\sqrt{\mathbf{1}^\top \mathbf{1}}=\sqrt{d}$.
    
    To show the containment, notice that by Lemma~\ref{lem:g_H_at_e}
    \begin{align*}
        \|x-e\|_e^2 = (M(x-e))^\top M(x-e) = \|Mx-\mathbf{1}\|_2^2 < 1,
    \end{align*}
    which implies that $(M x)_i > 0$ for all $i \in [d]$. Then, by Lemma~\ref{lem:lambda_l_intersect}, $x \in \Lambda_{++}$, and hence the proposition is proved.
\end{proof}

\subsection{Closed quadratic cone}
\label{sec:background_hp:closed_quadratic}

\begin{definition}[Closed quadratic cone]\label{def:quad_cone}
For any $e$ and for positive value $\alpha$, we define the closed quadratic cone as
\begin{align*}
    K_e(\alpha):=&~\{x\in \mathcal{E}:\ip{e}{x}_e\geq \alpha\|x\|_e\}\\
    = &~ \{x\in \mathcal{E}: \angle_e(e,x)\leq \arccos(\alpha / \sqrt{d})\},
\end{align*}
where $\angle_e$ denotes the angle between $e$ and $x$ measured in terms of the local inner product $\ip{\cdot}{\cdot}_e$.
\end{definition}

\begin{lemma} \label{lem:quad_cone_convex}
$K_e(\alpha)$ has the following properties:
\begin{itemize}
    \item $K_e(\alpha)$ is convex.
    \item The interior of $K_e(\alpha)$ contains the line segment (except the end points) of any two vectors $x, y$ in $K_e(\alpha)$, unless $x, x'$ are linearly dependent.
\end{itemize}

\end{lemma}
\begin{proof}
(1) Let $x,y\in K_e(\alpha)$. Then, for $t\in [0,1]$, 
\begin{align}
    \ip{e}{tx+(1-t)y}_e =&~ t\ip{e}{x}_e+(1-t)\ip{e}{y}_e \notag \\
    \geq &~ t\alpha\|x\|_e+(1-t)\alpha\|y\|_e \notag \\
    = &~ \alpha(\|tx\|_e+\|(1-t)y\|_e)\notag \\
    \geq &~ \alpha\|tx+(1-t)y\|_e, \label{eq:quad_cone_convex}
\end{align}
where the second step follows from definition of $K_e(\alpha)$, and the last step follows from the triangle inequality. 
Hence, $tx+(1-t)y\in K_e(\alpha)$, implying that $K_e(\alpha)$ is convex.

(2) For $x, y \in K_e(\alpha)$ and $0 < t < 1$, consider Eq. \eqref{eq:quad_cone_convex}. Notice that the equal sign in the last step holds if and only if $tx$ and $(1-t)y$ are linearly dependent, which is equivalent to $x$ and $y$ are linearly dependent, thus as long as $x, y$ are linearly independent,
\begin{align*}
    \langle e, tx + (1-t)y \rangle_e > \alpha \|tx + (1-t)y\|_e
\end{align*}
implying the line segment (except the end points) of $x, y$ is contained in the interior of $K_e(\alpha)$.
\end{proof}

\begin{proposition}[Proposition 8.9 in~\cite{rs14}]\label{prop:quad_cone_containing}
    We assume $e\in \Lambda_{++}$. We have 
    \begin{align*}
        K_e(\sqrt{d-1})\subseteq \Lambda_{+}\subseteq K_e(1).
    \end{align*}
\end{proposition}
For the completeness, we still provide a proof.
\begin{proof}
    We consider $x$ such that $0\ne x\in K_e(\sqrt{d-1})$, i.e, consider
    \begin{align*}
        \langle e, x\rangle_e \geq \sqrt{d-1}\|x\|_e.
    \end{align*}
    Define $\bar{x}:=\frac{\langle e, x\rangle_e}{\|x\|_e^2}x$, which is a positive scaled copy of $x$. Then, in order to show that $x\in \Lambda_+$, we can just show that for $\bar{x}$. We have
    \begin{align*}
        \|\bar{x}-e\|_e^2 = &~ \|\bar{x}\|_e^2 -2\langle e, \bar{x}\rangle_e + d\\
        = &~ \frac{\langle e, x\rangle_e^2}{\|x\|_e^2}-2\frac{\langle e, x\rangle_e^2}{\|x\|_e^2} + d\\
        \leq &~ -(d-1)+d \\
        = &~ 1,
    \end{align*}
    where the third step follows from $x \in K_e(\sqrt{d-1})$.
    
    By Proposition~\ref{prop:e_norm} 
    , $\bar{x}\in \Lambda_{+}$. Thus, $K_e(\sqrt{d-1})\subseteq \Lambda_+$.
    
    Now we show that $\Lambda_+\subseteq K_e(1)$. We first consider a fixed $x \in \Lambda_+$. We use $L$ to denote the subspace spanned by $x$ and $e$, and we consider $M$ to be described as in Lemma~\ref{lem:lax_dim_2}. By Lemma~\ref{lem:g_H_at_e},
    \begin{align*}
        \langle e, x\rangle_e = \mathbf{1}^\top Mx \quad \text{and} \quad \|x\|_e = \|Mx\|_2.
    \end{align*}
    By Lemma~\ref{lem:lambda_l_intersect}, $(Mx)_i\geq 0$ for all $i\in [d]$. Hence,
    \begin{align*}
        \mathbf{1}^\top Mx = \| Mx\|_1 \geq \|Mx\|_2, 
    \end{align*}
    which implies $x\in K_e(1)$.
\end{proof}

\begin{remark}
For any value of $r$ satisfying $r\leq \sqrt{d}$, the smallest cone that contains the ball $B_e(e,r)$ is $K_e(\sqrt{d-r^2})$.

For any $\alpha_1\leq \alpha_2$, we have $K_e(\alpha_1)\supseteq K_e(\alpha_2)$. Moreover, $\alpha=1$ is the largest value which for all hyperbolic polynomials of degree $d$, $K_e(\alpha)\supseteq \Lambda_+$. $\alpha=\sqrt{d-1}$ is the smallest value such that $K_e(\alpha)\subseteq \Lambda_+$. 
\end{remark}

\subsection{The dual of closed quadratic cones}\label{sec:background_hp:closed_quadratic_dual}

We also use the dual cone of $K_e(\alpha)$.

\begin{proposition}[Proposition 8.10 in~\cite{rs14}]\label{prop:dual_quad_cone}
    If $e\in \Lambda_{++}$, and $0<\alpha<\sqrt{d}$, then
    \begin{align*}
        K_e(\alpha)^*=\{H(e)s:s\in K_e(\sqrt{d-\alpha^2})\}.
    \end{align*}
\end{proposition}
For the completeness, we still provide a proof.
\begin{proof}
    According to Definition~\ref{def:quad_cone},
    \begin{align*}
        K_e(\alpha)=\{x:\angle_e(e,x)\leq \arccos(\alpha/\sqrt{d})\}.
    \end{align*}
    We have dual cone of the inner product $\ip{\cdot}{\cdot}_e$ to be
    \begin{align*}
        K_e(\alpha)^{*e} := &~ \{s:\ip{x}{s}_e\geq 0~\text{for all }x\in K_e(\alpha)\}\\
        = &~ \{s: \angle_e(e,s)\leq \arcsin(\sqrt{d}/\alpha)\}\\
        = &~ \{s: \angle_e(e,s)\leq \arccos(\sqrt{d-\alpha^2}/\sqrt{d})\}\\
        = &~ K_e(\sqrt{d-\alpha^2}),
    \end{align*}
    where the first step follows from the definition of dual cone, the second step follows from geometry, the third step follows from the property of inverse trigonometric function, and the last step follows from the definition of $K_e(\sqrt{d - \alpha^2})$.

    Since $K_e(\alpha)^*=\{s:\ip{x}{s}\geq 0~\text{for all }x\in K_e(\alpha)\}$ and $\ip{x}{s}_e=\ip{x}{H(e)s}$, we get that \begin{align*}
        K_e(\alpha)^* = &~ \{s : \langle x, s \rangle \ge 0~{\rm for~all}~ x \in K_e(\alpha) \} \\
        = &~ \{H(e)s: s\in K_e(\alpha)^{*e}\} \\
        = &~ \{H(e)s: s\in K_e(\sqrt{d - \alpha^2})\},
    \end{align*}
    where the first step follows from the definition of dual cone, second step follows from $\langle x, s \rangle_e = \langle x, H(e)s \rangle$, and the last step follows from $K_e(\alpha)^{*e} = K_e(\sqrt{d - \alpha^2})$.

    Hence completes the proof.
\end{proof}

\subsection{Eigenvalues of hyperbolic polynomials}\label{sec:background_hp:eigenvalues}

For a hyperbolic polynomial $p$, we can define the hyperbolic eigenvalue of an arbitrary vector $x \in \R^n$ similar with the eigenvalue of a matrix: the hyperbolic eigenvalues of $x$ are the positions of the zero points on the straight line passing through $x$ and with direction $e$.

\begin{definition}[Hyperbolic eigenvalue] 
Given a hyperbolic polynomial $p$ which is hyperbolic at the direction of $e \in \R^n$, for a vector $x$, we define the hyperbolic eigenvalues $\lambda_1(x) \ge \lambda_2(x) \ge \cdots \ge \lambda_d(x)$ of $x$ to be the roots of the function $p(te-x)$ as a function of $t$.

Usually we describe the hyperbolic eigenvalues of a vector $x$ as a single vector 
\begin{align*} 
\lambda = (\lambda_1(x)~\lambda_2(x)~\cdots~\lambda_d(x))^\top.
\end{align*}
\end{definition}

The following lemma shows that we can efficiently compute the first four moments of eigenvalues, which is very useful in the main algorithm.

\begin{lemma}[Evaluating the moments of $\lambda_e(x)$]\label{lem:eigen_moments}
We can evaluate the following moments at $\lambda_e(x)$ in $O(n^4)$ time for a degree-$d$ polynomial:
\begin{align*}
    \sum_{j\in [d]}\lambda_j,\quad \sum_{j\in [d]}\lambda_j^2,\quad \sum_{j\in [d]} \lambda_j^3 \quad \text{and} \quad \sum_{j\in [d]} \lambda_j^4. 
\end{align*}
Moreover, if we have an evaluation oracle for $p$ and take $\mathcal{T}_{O}$ time per oracle call, then we can compute these moments in $O(d^2+d\cdot\mathcal{T}_O)$ time.
\end{lemma}
\begin{proof}
We first argue for the $O(n^4)$ running time without an evaluation oracle, then we show how to further improve the efficiency using the oracle.

{\bf $O(n^4)$ method without an oracle.} 

We first compute the five leading coefficients $a_d,\dots, a_{d-4}$ of the univariate polynomial
\begin{align*}
    t\mapsto p(x+te)=\sum_{i=0}^d a_i t^i.
\end{align*}

We need to examine all of the monomials of degree at least $d-4$, which could be a large set. Actually, by the homogeneity of $p$, it is equivalent to compute the five trailing coefficients of $s\mapsto p(sx+e)$, since
\begin{align*}
    p(sx+e)
    = & ~ p(s(x+s^{-1}e)) \\
    = & ~ s^dp(x+s^{-1}e) \\
    = & ~ s^d\sum_{i=0}^d a_i s^{-i} \\
    = & ~ \sum_{i=0}^d a_i s^{d-i},
\end{align*}
where the second step follows from $p$ is homogeneous, and the third step follows from $p(x+te) = \sum_{i=0}^d a_i t^i$.

Hence, we can just compute the coefficients monomials of $p(sx+e)$ with degree at most 4. The total number of such monomials is at most
\begin{align*}
    \sum_{i=0}^4 \binom{n}{i}=\Theta(n^4).    
\end{align*}
Hence, we can compute $a_d,\dots, a_{d-4}$ in time $O(n^4)$.

{\bf $O(d^2+d\cdot {\cal T}_O)$ method with an oracle.} If we have oracle access to the polynomial $p$, then we can compute the coefficients by interpolation. More specifically, let $ \omega_d:=e^{2\pi \i / d}$ be the $d$-th root of unity. Then, we evaluate $p(x+te)$ at $\omega_d, \omega_d^2, \dots, \omega_d^{d}$, and get
\begin{align*}
    p(x+\omega_d^k e) - p(x) = \sum_{i=0}^d a_i \omega_d^{k(d-i)} - a_0 = \sum_{i=1}^d a_i \omega_d^{k(d-i)},
\end{align*}
thus
\begin{align*}
    \begin{bmatrix}
        \omega_d^d & \omega_d^{d-1} & \cdots & \omega_d\\
        \omega_d^{2\cdot d} & \omega_d^{2 \cdot(d-1) } & \cdots & \omega_d^2\\
        \vdots & \vdots & \ddots & \vdots\\
        \omega_d^{d\cdot d} & \omega_d^{d\cdot(d-1)} & \cdots & \omega_d^d
    \end{bmatrix}
    \begin{bmatrix}
    a_d\\
    a_{d-1}\\
    \vdots\\
    a_1
    \end{bmatrix}=
    \begin{bmatrix}
        p(x+\omega_d e)-p(x)\\
        p(x+\omega_d^2 e)-p(x)\\
        \vdots\\
        p(x+\omega_d^d e)-p(x)
    \end{bmatrix}
\end{align*}

We know that the inverse of Vandermonde matrix is
\begin{align}\label{eq:vandermonde}
    V:=\frac{1}{n}\begin{bmatrix}
    1 & 1 & \cdots & 1 & 1\\
    \omega_d & \omega_d^2 & \cdots & \omega_d^{d-1} & 1\\
    \vdots & \vdots & \ddots & \vdots & \vdots\\
    \omega_d^{d-1} & \omega_d^{2(d-1)} & \cdots & \omega_d^{(d-1)(d-1)} & 1
    \end{bmatrix}
\end{align}
And
\begin{align}\label{eq:interpolate}
    \begin{bmatrix}
    a_d\\
    \vdots\\
    a_1
    \end{bmatrix}= V\begin{bmatrix}
        p(x+\omega_d e)-p(x)\\
        \vdots\\
        p(x+\omega_d^d e)-p(x)
    \end{bmatrix}
\end{align}
Hence, if we can evaluate $p$ in $\mathcal{T}_O$ time, then we can compute $a_d, \dots, a_{d-4}$ in $O(d^2+d\cdot\mathcal{T}_O)$ time.

Then, by Newton-Girard identities,
\begin{align*}
    \sum_{j\in [d]}\lambda_j=&~ \frac{a_{d-1}}{a_d},\\
    \sum_{j\in [d]} \lambda_j^2=&~ \left(\frac{a_{d-1}}{a_d}\right)^2-\frac{a_{d-2}}{a_d},\\
    \sum_{j\in [d]} \lambda_j^3 = &~ \left(\frac{a_{d-1}}{a_d}\right)^3-\frac{3a_{d-1}a_{d-2}}{2a_d^2}+\frac{a_{d-3}}{2a_d},\\
    \sum_{j\in [d]} \lambda_j^4 = &~ \left(\frac{a_{d-1}}{a_d}\right)^4-\frac{2a_{d-1}^2a_{d-2}}{a_d^3}+\frac{a_{d-2}^2}{2a_d^2}+\frac{2a_{d-1}a_{d-3}}{3a_d^2}-\frac{a_{d-4}}{6a_d}.
\end{align*}
This step only takes $O(1)$ time, and hence, the lemma is proved. 

\end{proof}
\section{Tools for Primal Affine-Scaling Algorithm} \label{sec:app_affine_scaling}
In this section, we will introduce some results about the Primal Affine-Scaling Algorithm, and some tools for proving the main convergence theorem (Theorem~\ref{thm:rs_main}). 

We organize this section as follows. Section \ref{sec:QP_runtime} analyzes the running time to solve a quadratic program.
Sections \ref{sec:analysis_qt},  \ref{sec:analysis_wtqt} and \ref{sec:cone_geometry} give some useful technical lemmas/tools. Section \ref{sec:main_theorem_proof} proves the convergence theorem.

\subsection{Solving quadratic programs in cubic time 
} \label{sec:QP_runtime}

\begin{lemma}\label{lem:lp_time}
Given $A \in \R^{m \times n}$, $e \in {\rm Swath}(\alpha)$. If $g(e)$, $H(e)$ are known, then the optimal solution of a quadratic program $\QP_e(\alpha)$
can be found in time $O((n+m)^3)$.
\end{lemma}

\begin{proof}

Since $c$ 
is not in the row space of $A$, the optimal solution of the quadratic program is supposed to be at the boundary of $K_e(\alpha)$, i.e., 
\begin{align*}
    \ip{e}{x}_e^2-\alpha^2\|x\|_e^2 = 0.
\end{align*}
Then by Lemma \ref{lem:Hessian_at_e}, 
\begin{align*}
    \langle g(e),x\rangle^2-\alpha^2\langle x,H(e)x\rangle^2 = \ip{e}{x}_e^2-\alpha^2\|x\|_e^2 = 0.
\end{align*}

Now, by the first-order optimality conditions (Lagrange multiplier), for some vector $y \in \R^m$ and some scalar $\lambda \in \R$, 
the following equations are necessary for a point $x$ to be optimal for $\QP_e(\alpha)$:
\begin{align}\label{eq:QP_constraints}
    0 = & ~ \langle g(e),x\rangle^2-\alpha^2\langle x,H(e)x\rangle\notag\\
    Ax = & ~ b\\
    0 = & ~ \lambda c - A^\top y + \langle g(e),x\rangle g(e)-\alpha^2H(e)x \notag
\end{align}

Note that the second and third equations of Eq. \eqref{eq:QP_constraints} give a linear system in $(x,y,\lambda)$ with solutions which form a $1$-dimensional set (the second equation gives $m$ constraints and the third equation gives $n$ constraints, they are linearly independent because $A$ is full rank). For a linear system with $n+m$ equations and $n+m+1$ variables, it can be solved in $O((n+m)^3)$ by Gaussian elimination.  

By linearly parameterizing the $1$-dimensional set and substituting the variable $x$ which appeared at the 1st equation, we can get two candidates of solution, named $x',x''$, for which one of them is optimal for $\QP_e(\alpha)$. It only takes $O(n^2)$ time.

Finally, checking feasibility and objective values $\langle c, x'\rangle$ and $\langle c, x''\rangle$ reveals which is the optimal solution. 

Hence, the total running time is $O((n+m)^3)$.
\end{proof}

\begin{lemma}
Given $A \in \R^{m \times n}$ and $e \in {\rm Swath}(\alpha)$. Suppose $(x_e, y, \lambda)$ is the optimal solution to Eq. \eqref{eq:QP_constraints}, then the optimal solution $(y_e(\alpha), s_e(\alpha))$ of $\QP_e(\alpha)^*$ can be given in the following formula:
\begin{align}
    y_e(\alpha) = &~ -\frac{1}{\lambda}y,\label{eq:y_e}\\
    s_e(\alpha) = &~ \frac{\langle c, e-x_e\rangle}{d-\alpha^2}\left(-g(e)+\frac{\alpha^2}{\langle g(e),x_e\rangle}H(e)x_e\right)\notag\\
    = &~ \frac{\langle c, e-x_e\rangle}{d-\alpha^2}H(e)\left(e-\frac{\alpha^2}{\langle e,x_e\rangle_e}x_e\right)\label{eq:s_e}
\end{align}
\end{lemma}

\begin{proof}
Consider $(y_e,s_e)=(y_e(\alpha), s_e(\alpha))$ as the optimal solution of $\QP_e(\alpha)^*$. Let $(x_e,y,\lambda)$ be the solution of $\QP_e(\alpha)$. Then, by Eq.~\eqref{eq:QP_constraints},
\begin{align}\label{eq:first_order}
    A^\top y-\langle g(e),x_e\rangle g(e)+\alpha^2H(e)x_e=\lambda c.
\end{align}

Multiplying both sides by $(e-x_e)^\top$ gives
\begin{align*}
    (e-x_e)^\top A^\top y-\langle g(e),x_e\rangle \langle g(e), e-x_e\rangle + \alpha^2 (e-x_e)^\top H(e)x_e=\lambda\langle c, e-x_e\rangle.
\end{align*}

Since $Ax_e=b=Ae$ ($e\in \mathrm{Swath}(\alpha)$), $A(e-x_e) = 0$, we have
\begin{align*}
     -\langle g(e),x_e\rangle \langle g(e), e\rangle +\langle g(e),x_e\rangle^2 -\alpha^2\langle x_e, H(e)x_e\rangle+ \alpha^2 e^\top H(e)x_e = \lambda\langle c, e-x_e\rangle.
\end{align*}

By Eq.~\eqref{eq:QP_constraints}, $\langle g(e),x\rangle^2-\alpha^2\langle x,H(e)x\rangle=0$ and hence,
\begin{align*}
    -\langle g(e),x_e\rangle \langle g(e), e\rangle + \alpha^2 e^\top H(e)x_e = \lambda\langle c, e-x_e\rangle.
\end{align*}

By Lemma~\ref{lem:Hessian_at_e} and Proposition~ \ref{prop:e_norm}, we have
\begin{align*}
    -\langle g(e), e\rangle = \langle e, H(e)e\rangle = \|e\|_e^2=d,~\text{and}\quad e^\top H(e)x_e=-\langle g(e), x_e\rangle, 
\end{align*}
which gives
\begin{align}\label{eq:lambda}
    \lambda = \frac{(d-\alpha^2)\langle g(e), x_e\rangle}{\langle c, e-x_e\rangle}
\end{align}
Putting Eq.~\eqref{eq:lambda} into Eq.~\eqref{eq:first_order} 

Now we prove $y_e$ and $s_e$ are given in Eq. \eqref{eq:y_e} and \eqref{eq:s_e}:

First we show that $s_e(\alpha)\in K_e(\alpha)^*$. By Proposition~\ref{prop:dual_quad_cone}, we only need to show that 
\begin{align*}
    \frac{\langle c, e-x_e\rangle}{d-\alpha^2}\left(e-\frac{\alpha^2}{\langle e,x_e\rangle_e}x_e\right) \in K_e(\sqrt{d-\alpha^2}),
\end{align*}
which is equivalent to 
\begin{align*}
    \left\langle e, \frac{\langle c, e-x_e\rangle}{d-\alpha^2}\left(e-\frac{\alpha^2}{\langle e,x_e\rangle_e}x_e\right) \right\rangle_e \ge \sqrt{d - \alpha^2} \cdot \left\| \frac{\langle c, e-x_e\rangle}{d-\alpha^2}\left(e-\frac{\alpha^2}{\langle e,x_e\rangle_e}x_e\right) \right\|_e
\end{align*}

On the one hand,
\begin{align}
    \left\langle e, \frac{\langle c, e-x_e\rangle}{d-\alpha^2}\left(e-\frac{\alpha^2}{\langle e,x_e\rangle_e}x_e\right) \right\rangle_e=&~ \frac{\langle c, e-x_e\rangle}{d-\alpha^2}\langle e,e\rangle_e-\frac{\langle c, e-x_e\rangle}{d-\alpha^2}\frac{\alpha^2\langle e,x_e\rangle_e}{\langle e,x_e\rangle_e}\notag\\
    = &~ \frac{\langle c, e-x_e\rangle}{d-\alpha^2} d - \frac{\langle c, e-x_e\rangle}{d-\alpha^2}\alpha^2  \notag\\
    = &~ \langle c,e-x_e\rangle.\label{eq:e_se_ip}
\end{align}
where the second step follows from Proposition~\ref{prop:e_norm}. 

On the other hand,
\begin{align*}
    \left\langle e-\frac{\alpha^2}{\langle e,x_e\rangle_e}x_e, e-\frac{\alpha^2}{\langle e,x_e\rangle_e}x_e\right\rangle_e=&~ \langle e,e\rangle_e-2\frac{\alpha^2\langle e,x_e\rangle_e}{\langle e,x_e\rangle_e}+\frac{\alpha^4\langle x_e,x_e\rangle_e}{\langle g(e),x_e\rangle^2}\\
    = &~ d-2\alpha^2+\alpha^2\\
    =&~ d-\alpha^2,
\end{align*}
where the second step follows from Eq. \eqref{eq:QP_constraints}, $x_e \in \partial K_e(\alpha)$. Thus

\begin{align*}
    \sqrt{d-\alpha^2} \left\| \frac{\langle c, e-x_e\rangle}{d-\alpha^2}\left(e-\frac{\alpha^2}{\langle e,x_e\rangle_e}x_e\right) \right\|_e = &~ \frac{\langle c, e-x_e\rangle}{\sqrt{d-\alpha^2}} \cdot \sqrt{\left\langle e-\frac{\alpha^2}{\langle e,x_e\rangle_e}x_e, e-\frac{\alpha^2}{\langle e,x_e\rangle_e}x_e\right\rangle} \\
    = &~ \frac{\langle c, e-x_e\rangle}{\sqrt{d-\alpha^2}}\cdot \sqrt{d-\alpha^2} \\
    = &~ \langle c, e-x_e\rangle,
\end{align*}
which implies $\frac{\langle c, e-x_e\rangle}{d-\alpha^2}\left(e-\frac{\alpha^2}{\langle e,x_e\rangle_e}x_e\right) \in K_e(\sqrt{d-\alpha^2})$. Hence $s_e(\alpha) \in K_e(\alpha)^*$, and $(y_e,s_e)$ is a feasible solution for $\QP_e(\alpha)^*$.

Then, we show $(y_e, s_e)$ is optimal, since
\begin{align}
    \langle x_e, s_e\rangle = &~ \frac{\langle c, e-x_e\rangle}{d-\alpha^2}\left(\langle x_e, e\rangle_e-\frac{\alpha^2\langle x_e, x_e\rangle_e}{\langle e,x_e\rangle_e}\right)\notag\\
    = &~ \frac{\langle c, e-x_e\rangle}{d-\alpha^2}\left(\langle x_e, e\rangle_e-\langle e,x_e\rangle_e\right) \notag\\
    = &~ 0,\label{eq:x_s_complement}
\end{align}
where the second step follows from $x_e\in \partial K_e(\alpha)$.
And
\begin{align}\label{eq:strong_duality_QP}
    \langle b, y_e\rangle = &~ x_e^\top A^\top y_e\notag\\
    = &~ x_e^\top (c-s_e)\notag\\
    = &~ \langle c, x_e\rangle.
\end{align}
By strong duality, $x_e$ is optimal for $\QP_e(\alpha)$ and $(y_e, s_e)$ is optimal for $\QP_e(\alpha)^*$.

Thus, we complete the proof.
\end{proof}

\begin{corollary}
    Assume $0 < \alpha < 1$, $e \in {\rm Swath}_e(\alpha)$, then the duality gap between $\QP_e(\alpha)$ and $\QP_e(\alpha)^*$ is 
    \begin{align*}
        {\rm gap} = \langle c, e - x_e \rangle.
    \end{align*}
\end{corollary}
\begin{proof}
    Assuming $0<\alpha\leq 1$, $\Lambda_+\subseteq K_e(\alpha)$, and thus we have $K_e(\alpha)^*\subseteq \Lambda_+^*$. As a result, we have that, $(y_e(\alpha), s_e(\alpha))$ is feasible to the dual problem of the hyperbolic program $\HP$, which can be written as
    \begin{align*}
        \left.\begin{array}{cl}
        \max_y & \langle b, y\rangle\\
        \textrm{s.t.} & A^\top y+s=c\\
        & s\in \Lambda_+^*
        \end{array}\right\}\HP^*
    \end{align*}
    
    We know that $e$ is feasible for HP. The duality gap between $(y_e,s_e)$ and $e$ is
    \begin{align}\label{eq:gap_HP}
        \mathrm{gap}_e:=&~\langle c,e\rangle - \langle b,y_e\rangle \notag\\
        = &~\langle c, e- x_e\rangle,
    \end{align}
    where the last step follows from Eq. \eqref{eq:strong_duality_QP}. 
\end{proof}

\subsection{Analysis of a quadratic polynomial \texorpdfstring{$q(t)$}{}} \label{sec:analysis_qt}

In Section~\ref{sec:analysis_qt:quadratic_cone_matrix}, we present the quadratic cone of a matrix. In Section~\ref{sec:analysis_qt:moving_matrix}, we show how to move a matrix while keeping its properties.

\subsubsection{Quadratic cone of a matrix}\label{sec:analysis_qt:quadratic_cone_matrix}

\begin{definition}[Quadratic cone]
For matrix $E\succ 0$ and $\alpha>0$, the quadratic cone $K_E(\alpha)$ is defined as
\begin{align*}
    K_E(\alpha):=\{X\in \mathbb{S}^n: \tr [ E^{-1} X ]\geq \alpha \cdot ( \tr[ (E^{-1}X)^2 ] )^{1/2} \}.
\end{align*}
Hence, $X\in \partial K_I(\alpha)$ is equivalent to $\tr [ X ] =\alpha ( \tr [ X^2 ] ) ^{1/2}$. Note that for any $X \in K_{E}(\alpha)$, $X$ is a symmetric matrix. 
\end{definition}

\begin{proposition}[Lemma 3.1 in \cite{rs14}]\label{prop:dual_matrix_quad_cone}
$K_E(\alpha)^*=K_{E^{-1}}(\sqrt{d-\alpha^2})$. 
\end{proposition} 
For the completeness, we still provide a proof.
\begin{proof}
First consider the case when $E = I$. Notice that
\begin{align*}
    K_E(\alpha)^* = & ~ \{X : \forall~Y \in K_I(\alpha), \langle X, Y \rangle \ge 0 \} \\
    = & ~ \{X : \forall~Y~{\rm with}~\tr[Y] \ge \alpha\cdot \tr[Y^2]^{1/2}, \langle X, Y \rangle \ge 0 \} \\
    = & ~ \{X : \forall~Y~{\rm with}~\langle I, Y \rangle \ge \alpha \cdot \|Y\|_F , \langle X, Y \rangle \ge 0 \} \\
    = & ~ \{X : \forall~Y~{\rm with}~\angle(I, Y) \le \arccos{\frac{\alpha}{\sqrt{d}}}, \langle X,Y \rangle \ge 0 \} \\
    = & ~ \{X : \angle(I,X) \le \arcsin{\frac{\alpha}{\sqrt{d}}} \} \\
    = & ~ \{X : \angle(I,X) \le \arccos{\frac{\sqrt{d-\alpha^2}}{\sqrt{d}}} \} \\
    = & ~ \{X : \tr[X] \ge \sqrt{d-\alpha^2} \cdot \|X\|_F \} \\
    = & ~ K_I(\sqrt{d - \alpha^2}),
\end{align*}
where the first step follows from the duality, the second step follows from the definition of $K_I(\alpha)$, the third step follows from $\tr[Y] = \langle I, Y \rangle$ and $\tr[Y^2] = \|Y\|_F^2$, the fourth step follows from the law of cosines, the fifth step follows from geometry, the second step follows from the properties of inverse trigonometric function, the seventh step follows from the law of cosines, and the last step follows from the definition of $K_I(\sqrt{d-\alpha^2})$.

Then for general $E$. Notice that 
\begin{align}
    K_E(\alpha) = \{EX : X \in K_I(\alpha)\} \label{eq:quad_cone_affine}
\end{align}
thus 
\begin{align*}
    K_E(\alpha)^* = & ~ \{Y : \forall~X\in~K_E(\alpha), \langle X, Y \rangle \ge 0\} \\
    = & ~ \{Y : \forall~X\in~K_I(\alpha), \langle EX, Y \rangle \ge 0\} \\
    = & ~ \{Y : \forall~X\in~K_I(\alpha), \langle X, E^{-1}Y \rangle \ge 0\} \\
    = & ~ \{Y : E^{-1}Y \in K_I(\sqrt{d-\alpha^2}) \} \\
    = & ~ K_{E^{-1}}(\sqrt{d-\alpha^2}),
\end{align*}
where the first step follows from the duality, the second step follows from Eq. \eqref{eq:quad_cone_affine}, the third step follows from the properties of inner product, the fourth step follows from the definition of $K_I(\sqrt{d-\alpha^2})$, and the last step follows from Eq. \eqref{eq:quad_cone_affine}.

Hence completes the proof.

\end{proof}

\begin{proposition}[Lemma 4.4 in~\cite{rs14} for $E=I\in \R^{d\times d}$]\label{prop:trace_S}
For $X$ and $S$ defined in Lemma \ref{lem:quad_matrix}, we have
\begin{align*}
    \tr[ S ] = 1 \quad \text{and} \quad \tr[ X S ] = 0.
\end{align*}
\end{proposition}
\begin{proof}
By direct computation,
\begin{align*}
    \tr [ S ] =&~ \frac{1}{d-\alpha^2}(d-\alpha^2)= 1,
\end{align*}
and
\begin{align*}
    \tr [ X S ] =&~ \frac{1}{d-\alpha^2}\left( \tr [ X ] -\frac{\alpha^2 \tr[ X^2 ] }{\tr [ X ] }\right) = 0,
\end{align*}
where the second step follows from $X \in \partial K_I(\alpha)$.
\end{proof}

\begin{proposition}[Proposition 4.5 in~\cite{rs14} for $E=I\in \R^{d\times d}$]\label{prop:quad_condition}
    For $X$ and $S$ defined in Lemma \ref{lem:quad_matrix}, and for $t>-1$ and $0<\beta<1$, we have
    \begin{align*}
        E(t)\succ 0~\wedge S\in \mathrm{int}(K_{E(t)}(\beta)^*)\quad \Leftrightarrow \quad q(t)<\frac{1}{d-\beta^2}
    \end{align*}
    where we define $E(t)=\frac{1}{1+t}(I+tX)$. Further, $q$ is a strictly convex function. We have $q(0)=1/(d-\alpha^2)$ and $q'(0)<0$.
\end{proposition}
For the completeness, we still provide a proof.
\begin{proof}
By Proposition~\ref{prop:dual_matrix_quad_cone}, if $E(t)\succ 0$,
\begin{align*}
    S\in \mathrm{int}(K_{E(t)}(\beta)^*)=\mathrm{int}(K_{E(t)^{-1}}(\sqrt{d-\beta^2}))~\Leftrightarrow~\tr(E(t)S)> \sqrt{d-\beta^2}(\tr[(E(t)S)^2])^{1/2}
\end{align*}
For the left-hand side, using Proposition \ref{prop:trace_S},
\begin{align*}
    \tr [ E(t)S ] = & ~ \frac{1}{1 + t} ( \tr [ S ] + t \tr [X S] ) = \frac{1}{1+t}.
\end{align*}
For the right-hand side, according to the defition of $q(t)$,
\begin{align*}
    \left(\tr\left[(E(t)S)^2\right]\right)^{1/2}=\frac{1}{1+t}\left(\tr\left[((I+tX)S)^2\right]\right)^{1/2}=\frac{1}{1+t}\sqrt{q(t)}.
\end{align*}

Hence, if $E(t)\succ 0$, then 
\begin{align}\label{eq:q_s_first_half}
    S\in \mathrm{int}(K_{E(t)}(\beta)^*)\quad \Leftrightarrow &~ \tr[E(t)S] > \sqrt{d-\beta^2} (\tr[(E(t)S)^2])^{1/2} \notag \\
    \Leftrightarrow &~ \frac{1}{1+t} > \sqrt{d-\beta^2} \frac{\sqrt{q(t)}}{1+t} \notag \\
    \Leftrightarrow &~ q(t)<\frac{1}{d-\beta^2}.
\end{align}

The remaining is to show that when $q(t)<\frac{1}{d-\beta^2}$ and $t>-1$, $E(t)\succ 0$. Actually, we know that
\begin{align*}
    q(t)=\frac{\tr[ ( E( t ) S )^2 ] }{ ( \tr[ E(t) S ] )^2}=\frac{\sum_j \lambda_j^2}{(\sum_j \lambda_j)^{2}},
\end{align*}
where $\lambda_j$ are $S^{1/2}E(t)S^{1/2}$'s eigenvalues. For a fixed $S\succ 0$, if $E(t)$ has $k$ non-zero eigenvalues, by Cauchy–Schwarz inequality,
\begin{align*}
    q(t)=\frac{\sum_{j\in[d]} \lambda_j^2}{(\sum_{j\in [d]} \lambda_j)^{2}}\geq \frac{1}{k}.
\end{align*}
Consider the characteristic polynomial to $E(t)$:
\begin{align*}
    \det\left(\lambda I - \frac{1}{1+t}(I+tX)\right)\propto \det\left(\frac{\lambda(1+t)-1}{t}I-X\right)=\prod_{i=1}^d \left(\frac{\lambda(1+t)-1}{t}-\lambda_i\right) 
\end{align*}
Hence,
\begin{align*}
    \lambda(E(t))=\frac{1}{1+t}(\lambda(X)\cdot t + 1).
\end{align*}
Since $\tr[X]=\alpha (\tr[X^2])^{1/2}$, some eigenvalues of $X$ are positive, some are negative. Therefore, there exists $-1<\epsilon_- <0<\epsilon_+<1$, such that for $t\in (\epsilon_-, \epsilon_+)$, $E(t)\succ 0$. 
For $t=\epsilon_-$ or $t=\epsilon_+$, $E(t)\succeq 0$ with some zero eigenvalues. Hence, the number of non-zero eigenvalues $k\leq d-1$, and
\begin{align}\label{eq:q_0_1}
    q(t)\geq \frac{1}{d-1}~~~\text{for some }t\in (0,1).
\end{align}

Now, we need to consider the property of $q(t)$. First, we expand it in terms of $X$:

\begin{align*}
    q(t)=&~ \tr\left[((I+tX)S)^2\right]=\tr[(S+tXS)^2]=\tr[S^2]+2t\tr[XS^2]+t^2\tr[XSXS].
\end{align*}
For the first term,
\begin{align*}
    \tr[S^2]=&~ \frac{1}{(d-\alpha^2)^2}\left(d-2\alpha^2+\frac{\alpha^4\tr[X^2]}{(\tr[X])^2}\right)\\
    =&~ \frac{1}{(d-\alpha^2)^2}\left(d-2\alpha^2+\alpha^2\right)\\
    = &~ \frac{1}{d-\alpha^2},
\end{align*}
where the second step follows from $\alpha^2 \tr[X^2] = \tr[X]^2$.

For the second term,
\begin{align*}
    2\tr[XS^2]=&~ \frac{2}{(d-\alpha^2)^2}\tr\left[X-2\frac{\alpha^2X^2}{\tr[X]}+\frac{\alpha^4X^3}{(\tr[X])^2}\right]\\
    =&~ \frac{2}{(d-\alpha^2)^2}\left(-\tr[X]+\frac{\alpha^4\tr[X^3]}{\tr[X]^2}\right),
\end{align*}
where the second step follows from $\alpha^2 \tr[X^2] = \tr[X]^2$.

For the third term,
\begin{align*}
    \tr[XSXS]=&~ \frac{1}{(d-\alpha^2)^2}\tr\left[ X^2-2\frac{\alpha^2}{\tr[X]}X^3 + \frac{\alpha^4}{\tr[X]^2}X^4 \right]\\
    = &~ \frac{1}{(d-\alpha^2)^2}\left( \tr[X^2]- 2\frac{\alpha^2}{\tr[X]}\tr[X^3] +\frac{\alpha^4}{\tr[X]^2}\tr[X^4] \right).
\end{align*}

Hence,
\begin{align}\label{eq:q_expression}
    q(t) = \frac{1}{(d-\alpha^2)^2}\left( \tr[X^2]- 2\frac{\alpha^2}{\tr[X]}\tr[X^3] +\frac{\alpha^4}{\tr[X]^2}\tr[X^4] \right) t^2 \\
    + \frac{2}{(d-\alpha^2)^2}\left(-\tr[X]+\frac{\alpha^4\tr[X^3]}{(\tr[X])^2}\right) t +  \frac{1}{d-\alpha^2}\notag.
\end{align}
Therefore, we know that 
\begin{align}\label{eq:q_zero}
    q(0)=\frac{1}{d-\alpha^2}<\frac{1}{d-1}.
\end{align}
Further, by norm inequality, $\tr[X^3]\leq (\tr[X^2])^{3/2}=\alpha^{-3}(\tr[X])^3$. Thus, the linear term
\begin{align*}
    -\tr(X)+\frac{\alpha^4\tr[X^3]}{\tr[X]^2}\leq  -\tr[X] + \alpha\tr[X] < 0
\end{align*}
for $0<\alpha < 1$. And we get that
\begin{align}\label{eq:q_d_zero}
    q'(0) < 0.
\end{align}

By the property of quadratic polynomial, using Eq.~\eqref{eq:q_0_1}, Eq.~\eqref{eq:q_zero} and Eq.~\eqref{eq:q_d_zero}, we know that $q(t)$ is a strictly convex quadratic polynomial with a minimum point in $(0, \epsilon_+)$. Thus, when $q(t)<\frac{1}{d-1}$, we have $t\in (\epsilon_-, \epsilon_+)$, and $E(t)\succ 0$.
Therefore, combining with Eq.~\eqref{eq:q_s_first_half},
\begin{align*}
    E(t)\succ 0~\wedge S\in \mathrm{int}(K_{E(t)}(\beta)^*)\quad \Leftrightarrow \quad q(t)<\frac{1}{d-\beta^2}.
\end{align*}
And we complete the proof of the proposition.
\end{proof}

\subsubsection{Moving a matrix while keeping its properties}\label{sec:analysis_qt:moving_matrix}

\begin{lemma}[Proposition 10.3 in~\cite{rs14}]\label{lem:quad_matrix}
For a real number $0 < \alpha < 1$ and a $d \times d$ matrix $X \in \partial K_I (\alpha)$ and $X \ne {\bf 0}_{d \times d}$.
Let 
\begin{align*}
    S:= \frac{ 1 }{ d - \alpha^2 } \left( I - \frac{ \alpha^2 }{ \tr [ X ] } X \right).
\end{align*}
Define quadratic polynomial of $t$
\begin{align*}
    q(t) = \tr [ ( (I + t X) S )^2 ].
\end{align*}
Then $q(t)$ is strictly convex, and its minimizer $\bar{t}$ satisfies $\bar{t} > \frac{1}{2}\alpha / \|X\|_F$, and  if we define $\delta:=\bar{t}-\frac{1}{2}\alpha/\|X\|_F$, then for $t \in [\bar{t}-\delta, \bar{t}+\delta]$, we have
\begin{align*}
    E(t)\succ 0 ~\text{and}~S\in \mathrm{int}(K_{E(t)}(\beta)^*),
\end{align*}

where we define $E(t)=\frac{1}{1+t}(I+tX)$ and $\beta=\alpha\sqrt{(1+\alpha)/2}$.
\end{lemma}
For the completeness, we still provide a proof.

\begin{proof}[Proof]
By Eq.~\eqref{eq:q_expression}, the minimizer of $q(t)$ is 
\begin{align*}
    \bar{t}=&~ \frac{\frac{2}{(d - \alpha^2)^2} (-\tr[X] + \frac{\alpha^4 \tr[X^3]}{\tr[X]^2})} {-\frac{2}{(d - \alpha^2)^2} (\tr[X^2] - 2\frac{\alpha^2}{\tr[X]} \tr[X^3] + \frac{\alpha^4}{\tr[X]^2} \tr[X^4])}  \\ 
    =&~ \frac{(\tr[X])^3-\alpha^4\tr[X^3]}{\tr[X^2](\tr[X])^2-2\alpha^2\tr[X^3]\tr[X]+\alpha^4\tr[X^4]}\\
    \geq &~ \frac{\tr[X]^3-\alpha^4\tr[X^3]}{\tr[X^2](\tr[X])^2-2\alpha^2\tr[X^3]\tr[X]+\alpha^4(\tr[X^2])^2}\\
    = &~ \frac{(\tr[X])^3-\alpha^4\tr[X^3]}{(1+\alpha^{-2})(\tr[X])^4-2\alpha^2\tr[X^3]\tr[X]}\\
    = &~ \frac{\alpha^2}{2\tr[X]}\cdot \frac{\alpha^{-2}(\tr[X])^3-\alpha^2\tr[X^3]}{\frac{1+\alpha^{-2}}{2}(\tr[X])^3-\alpha^2\tr[X^3]}\\
    > &~ \frac{\alpha^2}{2\tr[X]}\\
    =&~ \frac{\alpha}{2(\tr[X^2])^{1/2}},
\end{align*}
where the third step follows from $\tr[X^4]\leq \tr[X^2]^2$ by Lemma \ref{lem:norm_inequality}
, and the sixth step follows from $\alpha^{-2}<1$ and $\tr[X]>0$. Therefore, the minimizer 
\begin{align*}
    \bar{t}>\frac{1}{2}\alpha / (\tr[X^2])^{1/2}=\frac{1}{2}\alpha / 2\|X\|_F.
\end{align*}

Let $\delta=\bar{t}-\frac{1}{2}\alpha/\|X\|_F$. By the property of quadratic polynomial, for $t\in [\bar{t}-\delta, \bar{t}+\delta]$, $q(t)\leq q(\frac{1}{2}\alpha/\|X\|_F)$. Now, we directly compute the value of $q(\frac{1}{2}\alpha/\|X\|_F)$. We first consider the linear and quadratic terms in Eq.~\eqref{eq:q_expression} (without the leading constant):

\begin{align*}
    &~ \left( \tr[X^2]- 2\frac{\alpha^2}{\tr[X]}\tr(X^3) +\frac{\alpha^4}{(\tr[X])^2}\tr[X^4] \right) \frac{\alpha^2}{4(\tr[X^2])} + 2\left(-\tr[X]+\frac{\alpha^4\tr[X^3]}{(\tr[X])^2}\right)\frac{\alpha}{2\tr[X^2]^{1/2}} \\
    = &~ \left( \tr[X^2]- 2\frac{\alpha^2}{\tr[X]}\tr(X^3) +\frac{\alpha^4}{(\tr[X])^2}\tr[X^4] \right)\frac{\alpha^4}{4(\tr[X])^2} + \left(-\tr[X]+\frac{\alpha^4\tr[X^3]}{(\tr[X])^2}\right)\frac{\alpha^2}{\tr[X]}\\
    = &~ \frac{\alpha^4}{4}-\alpha^2+\frac{\alpha^6\tr[X^3]}{2(\tr[X])^3}+\frac{\alpha^8\tr[X^4]}{4(\tr[X])^4}\\
    \leq &~ \frac{\alpha^4}{4}-\alpha^2+\frac{\alpha^3(\tr[X])^3}{2(\tr[X])^3}+\frac{\alpha^4(\tr[X])^4}{4(\tr[X])^4}\\
    = &~ \frac{\alpha^4}{2}+\frac{\alpha^3}{2}-\alpha^2,
\end{align*}
where the first step follows from $\tr[X]^2 = \alpha^2 \tr[X^2]$ and the third step follows from $\alpha^3 \tr[X^3] \le \tr[X]^3$ and $\alpha^4 \tr[X^4] \le \tr[X]^4$.

Hence,
\begin{align*}
    q\left(\frac{1}{2}\alpha/\|X\|_F\right) \le &~ \frac{1}{(d-\alpha^2)^2} (\frac{\alpha^4}{2} + \frac{\alpha^3}{2} - \alpha^2) + \frac{1}{d-\alpha^2} \\
    \leq &~ \frac{1}{d-\alpha^2}\left(1-\frac{2\alpha^2-\alpha^3-\alpha^4}{2(d-\alpha^2)}\right)\\
    < &~ \frac{1}{d-\alpha^2}\frac{1}{1+\frac{2\alpha^2-\alpha^3-\alpha^4}{2(d-\alpha^2)}}\\
    = &~ \frac{1}{d-\alpha^3(1+\alpha)/2} \\
    < &~ \frac{1}{d - \alpha^2(1+\alpha)/2},
\end{align*}
where the first step follows from Eq. \eqref{eq:q_expression}, the second step follows from $1-x < \frac{1}{1+x}$ for $x > 0$, and the last step follows from $0 < \alpha < 1$.

Define $\beta=\alpha\sqrt{\frac{1+\alpha}{2}}<1$, then 
\begin{align*}
    q\left(\frac{1}{2}\alpha/\|X\|_F\right)\leq &~ \frac{1}{d-\beta^2}.
\end{align*}

Therefore, for $t\in [\bar{t}-\delta, \bar{t}+\delta]$, $q(t)< \frac{1}{d-\beta^2}$. By Proposition~\ref{prop:quad_condition}, it implies $E(t)\succ 0$ and $S\in \mathrm{int}(K_{E(t)}(\beta)^*)$, which completes the proof.
\end{proof}

\subsection{Analysis of a quadratic polynomial \texorpdfstring{$\wt{q}(t)$}{}} \label{sec:analysis_wtqt}

In Section~\ref{sec:analysis_wtqt:relation}, we show the relationship between $\wt{q}(t)$ and $q(t)$. In Section~\ref{sec:analysis_wtqt:transformation}, we provide a transformation that turns vector to matrix. In Section~\ref{sec:analysis_wtqt:moving_vector}, we show how to move a vector while keeping its properties.

\subsubsection{The relation between \texorpdfstring{$\wt{q}(t)$}{} and \texorpdfstring{$q(t)$}{}}\label{sec:analysis_wtqt:relation}
 
\begin{lemma}[Proposition 10.4 in~\cite{rs14}]\label{lem:quad_poly_matrix}
Consider the polynomial $q$ in Lemma \ref{lem:quad_matrix}, and consider the polynomial $\wt{q}$ defined in Eq. \eqref{eq:quad_poly} which uses $\lambda(X)$, the eigenvalues of $X$, to act as $\lambda$. Then $q$ is a positive multiple of $\wt{q}$.
\end{lemma}
For the completeness, we still provide a proof.
\begin{proof}
    By Eq.~\eqref{eq:q_expression}, 
    \begin{align*}
        \frac{1}{(d-\alpha^2)^2}\cdot q(t) = \left( \tr[X^2]\tr[X]^2- 2\alpha^2\tr[X]\tr[X^3] +\alpha^4\tr[X^4] \right) t^2 \\
    + 2\left(-(\tr[X])^3+\alpha^4\tr[X^3]\right) t +  (d-\alpha^2)(\tr[X])^2.
    \end{align*}
    Writing trace as the sum of eigenvalues gives the same equation as Eq.~\eqref{eq:quad_poly}.
\end{proof}

\subsubsection{A transformation from vector to matrix}\label{sec:analysis_wtqt:transformation}

\begin{theorem}[Helton-Vinnikov Theorem~\cite{hv07,lpr05}]\label{thm:helton-vinnikov}
    Consider $e\in \Lambda_{++}$, and consider $L$ as a three-dimensional subspace of $\mathcal{E}$ for which contains the point $e$. We use $d$ to denote the degree of $p$. Then we have that, There is a linear transformation $T:L\rightarrow \mathbb{S}^d$ ($d\times d$ symmetric matrices) for which 
    \begin{align*}
        T(e)=I\quad \text{and}\quad p(x)=p(e)\det(T(x))~~\text{for all }x\in L.
    \end{align*}
\end{theorem}

\subsubsection{Moving a vector while keeping its properties}\label{sec:analysis_wtqt:moving_vector}

\begin{lemma}[Corollary 10.2 in~\cite{rs14}]\label{coro:quad_swath}
    Consider $e\in \mathrm{Swath}(\alpha)$ and $0<\alpha<1$. We use $\wt{q}$ to represent the quadratic polynomial defined as $\wt{q}(t)=at^2+bt+c$ as in Eq.~\eqref{eq:quad_poly}. Then we have, $\wt{q}$ is a strictly convex function. The minimizer of the quadratic polynomial $\wt{q}$, denoted as $t_e(\alpha)$, satisfies $\delta:=t_e(\alpha)-\frac{1}{2}\alpha/\|x_e(\alpha)\|_e>0$; besides, if
    \begin{align*}
        t_e(\alpha)-\delta \leq t\leq t_e(\alpha)+\delta\quad,
    \end{align*}
    then
    \begin{align*}
        \quad e(t)\in \mathrm{Swath}(\alpha)~\text{and}~s_e(\alpha)\in \mathrm{int}(K_{e(t)}(\beta)^*),
    \end{align*}
    
    where we have $\beta=\alpha\sqrt{(1+\alpha)/2}$ and $e(t)=\frac{1}{1+t}(e+tx_e(\alpha))$.
\end{lemma}

For the completeness, we still provide a proof.

\begin{proof}
Consider $e\in \Lambda_{++}$ and consider $L$ is a three-dimensional subspace of $\mathcal{E}$ such that it contains the point $e$. We define $T:L\rightarrow \mathbb{S}^d$ as a linear transformation as in Theorem~\ref{thm:helton-vinnikov}); hence we have $T(e)=I$.
    
First, observe that $p(\lambda e - x)=\det(T(\lambda e-x))=\det(\lambda I-X)$ where $X=T(x)$. Hence, $\lambda_e(x)=\lambda(X)$.

Then, consider the function $X\mapsto -\ln \det(X)$. By Lemma \ref{lem:mtx_det_partial}, since $X$ is symmetric,
the first derivative is
\begin{align*}
    \frac{\partial( -\ln \det(X))}{\partial X_{i,j}}=-(X^{-T})_{i,j} = -(X^{-1})_{i,j} ~~~ \forall (i,j)\in [d]\times [d].
\end{align*}
The second derivative is
\begin{align*}
    \frac{\partial^2( -\ln \det(X))}{\partial X_{i,j}\partial X_{k,\ell}}=-\frac{\partial (X^{-1})_{i,j}}{\partial X_{k,\ell}} = (X^{-1})_{i,k}(X^{-1})_{\ell, j}~~~\forall (i,j), (k,\ell)\in [d]\times [d].
\end{align*}
Hence, for any symmetric matrix $W\in \mathbb{S}^d$, the Hessian applied to $W$ is
\begin{align*}
    \nabla^2 (-\ln \det(X))[W]_{i,j}=\sum_{k,\ell} (X^{-1})_{i,k}(X^{-1})_{\ell ,j} W_{k,\ell} = (X^{-1}WX^{-1})_{i,j}~~~\forall (i,j)\in [d]\times [d].
\end{align*}
Note that, we use $H(X)$ to denote the Hessian of $-\ln \det(X)$. For $H(X)$ evaluate at $W$, we have 
\begin{align}\label{eq:hessian_det}
    H(X)[ W ] = X^{-1}WX^{-1}.
\end{align}

Therefore, for all $\bar{e}\in \Lambda_{++}\cap L$ and $u,v\in L$, 
\begin{align}
    \langle u,v\rangle_{\bar{e}} = &~ \langle u, \nabla^2(-\ln p(\bar{e}))v\rangle\notag \\
    = &~ \langle T(u), \nabla^2(-\ln \det(T(\bar{e}))) T(v)\rangle\notag \\
    = &~ \langle T(u), H(T(\bar{e}))T(v)\rangle\notag \\
    = &~ \tr[T(u)T(\bar{e})^{-1}T(v)T(\bar{e})^{-1}]\notag\\
    = &~ \tr[T(\bar{e})^{-1}T(u)T(\bar{e})^{-1}T(v)],\label{eq:norm_poly_mat}
\end{align}
where the second step follows from the definition of $H(e)$, the third step follows from  the forth step follows from Eq.~\eqref{eq:hessian_det}. 

A direct consequence is that, for $x\in L$ and $\bar{e}\in \Lambda_{++}\cap L$, and any $\gamma>0$,
\begin{align}\label{eq:cone_poly_mat}
    x\in K_{\bar{e}}(\gamma)\quad \Leftrightarrow \quad T(x)\in K_{\bar{E}}(\gamma).
\end{align}

Now we consider $x$ such that $0\ne x\in \partial K_e(\alpha)$. Thus we have $X=T(x)\in \partial K_I(\alpha)$. We use $\wt{q}$ to denote the quadratic polynomial described in Eq.~\eqref{eq:quad_poly}. By $\lambda_e(x)=\lambda(X)$, we have that, $\wt{q}$ is identical for the polynomial which is described as Lemma~\ref{lem:quad_poly_matrix}. Hence by the proposition, we know $\wt{q}$ is a positive multiplied copy of the polynomial $q$ described as Lemma~\ref{lem:quad_matrix}.

Hence, by Lemma~\ref{lem:quad_matrix}, $\wt{q}$ is a strictly convex function, whose minimizer $t_e(\alpha)$ satisfies the following
\begin{align*}
    \delta:=t_e(\alpha)-\frac{1}{2}\alpha / \|x\|_e = t_e(\alpha)-\frac{1}{2}\alpha / \|X\|_F>0.
\end{align*}
Besides, by $\Lambda_{++}\cap L = \{\bar{e}\in L: T(\bar{e})\succ 0\}$, from Lemma~\ref{lem:quad_matrix} we have
\begin{align*}
    t_e-\delta \leq t \leq t_e + \delta \quad \Rightarrow \quad e(t)\in \Lambda_{++}
\end{align*}
for which $e(t)= \frac{1}{1+t}(e+tx)$.

Now we consider that $e\in \mathrm{Swath}(\alpha)$ and $0<\alpha <1$. Let $x=x_e(\alpha)$ to be the optimal solution with respect to $\QP_e(\alpha)$. To establish Lemma~\ref{coro:quad_swath}, now we only need to show for $e(t)=\frac{1}{1+t}(e+tx_e)$ and $\beta=\alpha\sqrt{(1+\alpha)/2}$, 
\begin{align*}
    t_e-\delta \leq t \leq t_e +\delta \quad \Rightarrow \quad e(t)\in \mathrm{Swath}(\beta)\text{ and }s_e\in \mathrm{int}(K_{e(t)}(\beta)^*),
\end{align*}
where $(y_e(\alpha), s_e(\alpha))$ is the dual optimal solution. 

First, it's easy to show that $e(t)\in \mathrm{Swath}(\beta)$ since
\begin{align*}
    Ae(t)=\frac{1}{1+t}(Ae+tAx_e)=b.
\end{align*}

Also, $(y_e,s_e)$ is a strictly feasible solution for $\QP_{e(t)}(\beta)^*$ since $s_e\in \mathrm{int}(K_{e(t)}(\beta)^*)$. $e(t)$ is a strictly feasible solution for $\QP_{e(t)}(\beta)$ as well. By the strong duality theorem of conic programming, since the primal and dual problems are both strictly feasible, the optimal solutions are the same and attainable. Therefore, there is an optimal solution to $\QP_{e(t)}(\beta)$, that is, $e(t)\in \mathrm{Swath}(\beta)$.

Second, we need to show that, for each $t\in [t_e-\delta, t_e+\delta]$, 
\begin{align*}
    &s_e=\frac{\langle c, e-x_e\rangle}{d-\alpha^2}H(e)\left(e-\frac{\alpha^2}{\langle e,x_e\rangle_e}x_e\right)\text{ satisfies}\\
    &\langle z,s_e\rangle > 0 \text{ for all }0\ne z\in K_{e(t)}(\beta),
\end{align*}
where $s_e$ is defined in Eq. \eqref{eq:s_e}.

Since we know that $K_{e(t)}(\beta)^*=\{H(e(t))s:  \langle z,s\rangle_{e(t)} \ge 0, \forall z\in K_{e(t)}(\beta)\}$, equivalently, 
\begin{align*}
    &s:= \frac{1}{\langle c,e-x_e\rangle}H(e(t))^{-1}s_e=\frac{1}{d-\alpha^2}H(e(t))^{-1}H(e)\left(e-\frac{\alpha^2}{\langle e,x_e\rangle_e}x_e\right)~\text{satisfies}\\
    &\langle z,s\rangle_{e(t)} > 0 \text{ for all }0\ne z\in K_{e(t)}(\beta)
\end{align*}

For some fixed $t \in [t_e-\delta, t_e+\delta]$ and fixed $z$ such that $0\ne z\in K_{e(t)}(\beta)$. Consider $L$ as the three-dimensional subspace which is spanned by three vectors $e, x_e, z$. For $T$ described as Theorem~\ref{thm:helton-vinnikov}, assume $E(t) = T(e(t)), Z=T(z)$, and $X_e= T(x_e)$.
By Eq.~\eqref{eq:cone_poly_mat}, $0\ne Z\in K_{E(t)}(\beta)$. Also, $s\in L$, 
\begin{align*}
    S=&~ \frac{1}{d-\alpha^2}H(E(t))^{-1}H(I)\left(I-\frac{\alpha^2}{ \tr[X_e]}X_e\right)\\
    = &~ \frac{1}{d-\alpha^2} E(t)I^{-1}\left(I-\frac{\alpha^2}{ \tr[X_e]}X_e\right) I^{-1} E(t)\\
    = &~ \frac{1}{d-\alpha^2} E(t)\left(I-\frac{\alpha^2}{ \tr[X_e]}X_e\right) E(t)
\end{align*}
where the first step follows from $\langle e, x_e\rangle_e = \tr[X_e]$ by Eq.~\eqref{eq:norm_poly_mat}, and the second step follows from Eq.~\eqref{eq:hessian_det}. 
By Lemma~\ref{lem:quad_matrix},
\begin{align*}
    E(t)^{-1}SE(t)^{-1} = \frac{1}{d-\alpha^2} \left(I-\frac{\alpha^2}{ \tr[X_e]}X_e\right) \in \mathrm{int}(K_{E(t)}(\beta)^*).
\end{align*}
 Hence, $S\in \mathrm{int}(K_{E(t)}(\beta)^{*E(t)})=\{S: \langle S, Z\rangle_{E(t)}>0~\forall Z\in K_{E(t)}(\beta)\}$. Then, we have
\begin{align*}
    \langle z,s\rangle_{e(t)}=\tr[E(t)^{-1} Z E(t)^{-1}S]=\langle Z, S\rangle_{E(t)}>0.
\end{align*}
Thereby, $s_e\in \mathrm{int}(K_{e(t)}(\beta)^*)$ and hence completing the proof.
\end{proof}

\subsection{A useful conclusion of Cone geometry} \label{sec:cone_geometry}

Last thing we need in order to prove Theorem~\ref{thm:rs_main} is a proposition about cone geometry. 
\begin{proposition}[Proposition 5.1 in~\cite{rs14}]\label{prop:cone_geom}
Consider $e$ such that $\|e\|_e= \sqrt{d}$, and $0<\alpha < 1$. Let $\beta = \alpha \sqrt{\frac{1+\alpha}{2}}$. Assume $\bar{s}\in \mathrm{int}(K_e(\beta)^*)$, and assume $\bar{x}$ is the optimal solution of the optimization problem
\begin{align*}
    \begin{array}{cl}
    \min_x & \langle \bar{s}, x\rangle\\
    \mathrm{s.t.} & x\in e+L\\
    & x\in K_{e}(\alpha),
\end{array}
\end{align*}
such that $\|x\|_e\geq \sqrt{\frac{2d}{1-\alpha}}$. Then we have, there is $\bar{s}'\in (\bar{s}+L^\bot)\cap K_e(\alpha)^*$ such that
\begin{align*}
    \langle e, \bar{s}-\bar{s}'\rangle_e\geq \kappa \|\bar{s}\|_e,
\end{align*}
where $\kappa = \alpha \sqrt{\frac{1-\alpha}{8}}$.
\end{proposition}

\subsection{Proof of the main theorem} \label{sec:main_theorem_proof}
Here in this section, we provide the following detailed proof. 

\begin{proof}[Proof of Theorem~\ref{thm:rs_main}]
First, for $e\in \mathrm{Swath}(\alpha)$ and $e'\in \frac{1}{1+t_e}(e+t_e x_e)$, by Lemma~\ref{coro:quad_swath}, $e'\in \mathrm{Swath}(\beta)\subseteq \mathrm{Swath}(\alpha)$, since $\beta < \alpha$.

Second, $\langle c,e'\rangle < \langle c,e\rangle$ since $\langle c,x_e\rangle < \langle c, e\rangle$ and $e'$ is a linear combination of $x_e$ and $e$.

Third, by Lemma~\ref{coro:quad_swath}, we know that $s_e\in \mathrm{int}(K_{e'}(\alpha)^*)$, which implies $(y_e, s_e)$ is a feasible solution for the quadratic program $\QP_{e'}(\alpha)^*$. Hence, $\langle b, y_{e'}\rangle >\langle b, y_e\rangle$ by the duality. 

Last, we need to show that, for $j=0,1,\dots$,
\begin{align*}
    \frac{\mathrm{gap}_{e_{j+1}}}{\mathrm{gap}_{e_j}} & \leq 1-\frac{\kappa}{\kappa+\sqrt{d}}\\
     \text{for }j=i & \text{ or }j=i+1\text{ (possibly both)},
\end{align*}
where the parameter $\kappa = \alpha \sqrt{\frac{1-\alpha}{8}}$. Since $e_0$ is an arbitrary point in $\mathrm{Swath}(\alpha)$, it's enough to prove it for $j=0$. For brevity, let $x_i := x_{e_i}$, $s_i := s_{e_i}$, $y_i:=y_{e_i}$.

Note that
\begin{align*}
    \mathrm{gap}_{1}=\mathrm{gap}_0 - \langle c, e_0 - e_{1}\rangle - b^\top (y_{1}-y_0).
\end{align*}

Hence, we simply have
\begin{align*}
    \frac{\mathrm{gap}_{1}}{\mathrm{gap}_0}\leq &~1-\frac{\langle c, e_0 - e_{1}\rangle}{\mathrm{gap}_0},~~~\text{and}\\
    \frac{\mathrm{gap}_{1}}{\mathrm{gap}_0}\leq &~ 1-\frac{b^\top (y_{1}-y_0)}{\mathrm{gap}_0}.
\end{align*}
Let $t_0,t_1$ denote the minimizer of $\wt{q}_0$ and $\wt{q}_1$. Now, consider two cases:
\begin{itemize}
    \item Either $t_0>\kappa / \sqrt{d}$ or $t_1>\kappa/\sqrt{d}$.
    \item Both $t_0,t_1\leq \kappa / \sqrt{d}$.
\end{itemize}

In the first case, we can show that
\begin{align*}
    \frac{\mathrm{gap}_{1}}{\mathrm{gap}_0}\leq &~1-\frac{\langle c, e_0 - e_{1}\rangle}{\mathrm{gap}_0}\\
    = &~ 1-\frac{\langle c, e_0 - e_{1}\rangle}{\langle c, e_0-x_0\rangle }\\
    = &~ 1-\frac{\langle c, e_0-\frac{1}{1+t_0}(e_0+t_0x_0)\rangle}{\langle c, e_0-x_0\rangle}\\
    = &~ 1-\frac{t_0}{1+t_0}\\
    < &~ \frac{\sqrt{d}}{\kappa+\sqrt{d}},
\end{align*}
where the second step follows from Eq.~\eqref{eq:gap_HP}, and the last step follows from $t_0>\kappa/\sqrt{d}$.

In the second case, we use $L$ to denote the nullspace of $A$, and we use $L^\bot$ to denote the orthogonal complement of $L$. Also, assume $L^{\bot_{e_1}}$ as the orthogonal complement with respect to the inner product $\langle ~,~\rangle_{e_1}$, i.e.,
\begin{align*}
    L^{\bot_{e_1}} := &~ \{x: \langle x, v\rangle_{e_1}=0~\forall v\in L\}\\
    = &~ \{H(e_1)^{-1}x: x\in L^\bot\}.
\end{align*}
Since $e_1\in \mathrm{Swath}$, we have $Ae_1=b$ and we can write $\QP_{e_1}(\alpha)$ as
\begin{align*}
    \left.\begin{array}{cl}
    \min_x & \langle c, x\rangle\\
    \textrm{s.t.} & x\in e_1+L\\
    & x\in K_{e_1}(\alpha)
    \end{array}\right\} \QP_{e_1}(\alpha)
\end{align*}
Notice that $A^\top y_0 + s_0 = c$, and then $c-s_0\in L^\bot$
. Hence, for all feasible $x,x'$, $\langle c-s_0, x-x'\rangle=0$. Thus, $x_1$ is also optimal for
\begin{align*}
\begin{array}{cl}
    \min_x & \langle s_0, x\rangle\\
    \textrm{s.t.} & x\in e_1+L\\
    & x\in K_{e_1}(\alpha)
\end{array}
\end{align*}
Let $\bar{s}_0=H(e_1)^{-1}s_0$. Then, $\langle s_0, x\rangle = \langle \bar{s}_0, x\rangle_{e_1}$. Consequently, we can also express the above optimization problem as
\begin{align*}
\begin{array}{cl}
        \min_x & \langle \bar{s}_0, x\rangle_{e_1}\\
    \textrm{s.t.} & x\in e_1+L\\
    & x\in K_{e_1}(\alpha)
\end{array}
\end{align*}
By Lemma~\ref{coro:quad_swath}, $s_0\in \mathrm{int}(K_{e_1}(\beta))$. It follows that $\bar{s}_0\in \mathrm{int}(K_{e_1}(\beta)^{*{e_1}})$.

Since $t_1> \frac{1}{2}\alpha\|x_1\|_{e_1}$ and $t_1\leq \kappa/\sqrt{d}$, we have $\|x_1\|_{e_1}\geq \sqrt{\frac{2d}{1-\alpha}}$.
Now, by Proposition~\ref{prop:cone_geom}, there is a point $\bar{s}'\in (\bar{s}_0+L^{\bot_{e_1}})\cap K_{e_1}(\alpha)^{*e_1}$ such that
\begin{align}\label{eq:cone_geometry}
    \langle e_1, \bar{s}_0-\bar{s}'\rangle_{e_1}\geq \kappa \|\bar{s}_0\|_{e_1}.
\end{align}
Assume $s'=H(e_1)\bar{s}'$. By $\bar{s}'\in \bar{s}_0+L^{\bot_{e_1}}$, $s'\in s_0+L^\bot$, and hence $s'\in c+L^\bot$. Thus, there exists a $y'$ such that $A^\top y'+s'=c$. Moreover, from $\bar{s}'\in K_{e_1}(\alpha)^{*e_1}$ follows $s'\in K_{e_1}(\alpha)^*$, and thus $(y', s')$ is feasible for $\QP_{e_1}(\alpha)^*$.

Consequently,
\begin{align*}
    b^\top y_1-b^\top y_0 \geq &~ b^\top y'-b^\top y_0\\
    = &~ \langle Ae_1, y'-y_0\rangle\\
    = &~ \langle e_1, A^\top (y'-y_0)\rangle\\
    = &~ \langle e_1, s_0-s'\rangle\\
    = &~ \langle e_1, \bar{s}_0-\bar{s}'\rangle_{e_1}\\
    \geq &~ \kappa\|\bar{s}_0\|_{e_1},
\end{align*}
where the first step follows from $y_1$ is the optimal solution of $QP_{e_1}(\alpha)$, the second step follows from $b = Ae_1$, the third step follows from $\langle A x, y \rangle = y^\top A x = \langle x,  A^\top y \rangle$, 
 the fourth step follows from $A^\top y_0 + s_0 = c$, the fifth step follows from the definition of $\bar{s}_0$, and the last step follows from Eq. \eqref{eq:cone_geometry}.
and
\begin{align*}
     \|\bar{s}_0\|_{e_1} 
    \geq &~ \frac{1}{\sqrt{d}}\langle e_1, \bar{s}_0\rangle_{e_1} \\ 
    = &~ \frac{1}{\sqrt{d}}\langle e_1, s_0\rangle\\
    = &~ \frac{1}{\sqrt{d}}\frac{1}{1+t_0}\langle e_0+t_0x_0, s_0\rangle\\
    = &~ \frac{1}{\sqrt{d}(1+t_0)}\langle e_0, s_0\rangle \\ 
    = &~ \frac{1}{\sqrt{d}(1+t_0)}\langle c, e_0-x_0\rangle \\ 
    = &~ \frac{1}{\sqrt{d}(1+t_0)} \cdot \mathrm{gap}_0\\
    \geq &~ \frac{1}{\sqrt{d}(1+\kappa/\sqrt{d})}\cdot \mathrm{gap}_0 \\ 
    = &~ \frac{1}{\kappa +\sqrt{d}}\cdot \mathrm{gap}_0,
\end{align*}
where the first step follows from Cauchy-Schwarz inequality, the second step follows from the definition of $\bar{s}_0$, the third step follows from $e_1 = e_0 + t_0x_0$, the fourth step follows from Eq. \eqref{eq:x_s_complement}, the fifth step follows from Eq. \eqref{eq:e_se_ip}, the sixth step follows from the definition of gap, the seventh step follows from $t_0 \le \kappa / \sqrt{d}$, and the last step follows from simplification.

Thus 
\begin{align*}
    b^\top y_1-b^\top y_0 \geq &~ \frac{\kappa}{\kappa +\sqrt{d}}\cdot \mathrm{gap}_0.
\end{align*}
Therefore,
\begin{align*}
    \frac{\mathrm{gap}_{1}}{\mathrm{gap}_0}\leq  1-\frac{b^\top (y_{1}-y_0)}{\mathrm{gap}_0} \leq 1-\frac{\kappa}{\kappa + \sqrt{d}},
\end{align*}
which completes the proof of the theorem.
\end{proof}

\section{Gradient and Hessian}
\label{sec:grad_ad_hes}

In this section, we provide missing proofs for Section~\ref{sec:main_alg}, which concerns the number of iterations and the computation of gradient and Hessian. In Section~\ref{sec:grad_ad_hes:lem_iterations}, we provide proof of Lemma~\ref{lem:iterations}. In Section~\ref{sec:grad_ad_hes:lem_fast_grad}, we provide proof of Lemma~\ref{lem:fast_grad}. In Section~\ref{sec:grad_ad_hes:lem_fast_hessian}, we provide proof of Lemma~\ref{lem:fast_hessian}.

\subsection{Proof of Lemma~\ref{lem:iterations}}\label{sec:grad_ad_hes:lem_iterations}

\begin{lemma}[Restatement of Lemma~\ref{lem:iterations}]\label{lem:iterations_formal}
Assume that for the starting point $e_0$, $\mathrm{gap}_0\leq 1$. Then, for any $0<\delta<1$, Algorithm \ref{alg:main} uses $O( \sqrt{d} \log( 1 / \delta ) )$ iterations to make the duality gap at most $\delta$.
\end{lemma}
\begin{proof}
Suppose that after $t$ iterations, $\mathrm{gap}_t\leq \delta$. Then, we have
\begin{align*}
    \mathrm{gap}_t \leq~ \frac{\mathrm{gap}_t}{\mathrm{gap}_0} =~ \prod_{i=0}^{t-1}\frac{\mathrm{gap}_{i+1}}{\mathrm{gap}_i}.
\end{align*}
By Theorem~\ref{thm:rs_main}, $\frac{\mathrm{gap}_{i+1}}{\mathrm{gap}_i}\leq 1$ for all $i\in [t-1]$.  
Moreover, for all $i\in [t-2]$, either $\frac{\mathrm{gap}_{i+1}}{\mathrm{gap}_i}\leq 1-\frac{\kappa}{\kappa + \sqrt{d}}$ or $\frac{\mathrm{gap}_{i+2}}{\mathrm{gap}_{i+1}}\leq 1-\frac{\kappa}{\kappa + \sqrt{d}}$, where $0<\kappa<1$ is a constant. Hence, we have
\begin{align*}
    \prod_{i=0}^{t-1}\frac{\mathrm{gap}_{i+1}}{\mathrm{gap}_i} \leq &~ \left(1-\frac{\kappa}{\kappa+\sqrt{d}}\right)^{t/2}\\
    = &~ \left(1-\frac{\kappa}{\kappa+\sqrt{d}}\right)^{\frac{\kappa+\sqrt{d}}{\kappa}\frac{t\kappa}{2(\kappa+\sqrt{d})}}\\
    \leq &~ \exp\left(-\frac{t\kappa}{2(\kappa+\sqrt{d})}\right),
\end{align*}
where the last step follows from $(1 - 1/x)^x \le 1/e$ for $x > 1$.

If $t\geq \frac{2(\kappa+\sqrt{d})}{\kappa} \ln( 1 / \delta )$, then we have
\begin{align*}
   \mathrm{gap}_t \leq~   \prod_{i=0}^{t-1}\frac{\mathrm{gap}_{i+1}}{\mathrm{gap}_i} \leq ~ \delta.
\end{align*}
Therefore, after $O(\sqrt{d} \log ( 1 / \delta ) )$ iterations in Algorithm~\ref{alg:main}, the duality gap is at most $\delta$.
\end{proof}

\subsection{Proof of Lemma~\ref{lem:fast_grad}}\label{sec:grad_ad_hes:lem_fast_grad}

\begin{lemma}[Restatement of Lemma~\ref{lem:fast_grad}]\label{lem:fast_grad_formal}
For hyperbolic polynomial $p$ of degree $d$, given an evaluation oracle for $p$, for any $x, w\in \R^n$, $\langle \nabla (-\ln p(x)), w\rangle$ can be computed in $O(d\mathcal{T}_O)$ time, where $\mathcal{T}_O$ is the time per oracle call.
\end{lemma}
\begin{proof}
First, we have $\nabla(-\ln p(x))=\frac{-\nabla p(x)}{p(x)}$. Then, we just need to compute 
\begin{align*}
    \langle \nabla p(x), w\rangle = \frac{\d p(x+tw)}{\d t}\Big |_{t=0}.
\end{align*}
Note that we can write $p(x+tw)$ as $\sum_{i=0}^d a_i t^i$. Then, we have $\langle \nabla p(x), w\rangle = a_1$.

Next, we can use a similar method in Lemma~\ref{lem:eigen_moments} to compute the first coefficient of $p(x+tw)$ by calling the oracle $d$ times. More specifically, by Eq.~\eqref{eq:vandermonde} \eqref{eq:interpolate}, we have
\begin{align}\label{eq:dir_grad}
    a_1 = &~ \frac{1}{n}\sum_{j=1}^d \omega_d^{-j} \left(p(x + \omega_d^j w) - p(x)\right)\notag\\
    = &~ \frac{1}{n}\sum_{j=1}^d \omega_d^{-j} p(x + \omega_d^j w),
\end{align}
where the third line follows from $\sum_{j=1}^n \omega_d^{-j} = 0$. Then,
\begin{align*}
    \langle \nabla (-\ln p(x)), w\rangle=\frac{-a_1}{p(x)} = \frac{-1}{np(x)}\sum_{j=1}^d \omega_d^{-j} p(x + \omega_d^j w)
\end{align*}

Hence, $\langle \nabla (-\ln p(x)), w\rangle$ can be computed in $O(d\mathcal{T}_O)$ time.
\end{proof}

\subsection{Proof of Lemma~\ref{lem:fast_hessian}}\label{sec:grad_ad_hes:lem_fast_hessian}

\begin{lemma}[Restatement of Lemma~\ref{lem:fast_hessian}]\label{lem:fast_hessian_formal}
For hyperbolic polynomial $p$ of degree $d$, given an evaluation oracle for $p$, for any $w\in \R^n$, $\nabla^2 (-\ln p(x))w$ can be computed in $O(nd^2\mathcal{T}_O)$ time.
\end{lemma}

\begin{proof}
By simple calculation, for every $i\in [n]$, 
\begin{align*}
    \left(\nabla^2 (-\ln p(x))w\right)_i = & ~ \langle - \nabla \frac{\partial_i p(x)}{p(x)} , w \rangle\\
    = & ~ \frac{\partial_i p(x) \cdot \langle \nabla p(x), w\rangle}{p(x)^2}-\frac{\langle \nabla \partial_i p(x), w\rangle}{p(x)},
\end{align*} 
where $\partial_i p(x) = \frac{\partial p(x)}{\partial x_i}=\langle \nabla p(x), e_i\rangle$ for $e_i=\begin{bmatrix}0 &\cdots & 0 & 1& 0 & \cdots & 0\end{bmatrix}^\top$ the $i^{th}$ basis vector.

Then, for the first term, by Lemma~\ref{lem:fast_grad}, $\langle \nabla  p(x), e_i\rangle$ and $\langle \nabla p(x), w\rangle$ can be computed by calling the evaluation oracle $O(d)$ times.

Hence, this term can be computed in $O(d\mathcal{T}_O)$ time.

For the second term, by Eq.~\eqref{eq:dir_grad}, 
\begin{align*}
    \partial_i p(x) = \frac{\d p(x+te_i)}{\d t}\Big |_{t=0} 
    = \frac{1}{n}\sum_{j=1}^d \omega_d^{-j} p(x + \omega_d^j e_i).
\end{align*}
 Hence,
\begin{align*}
    \langle \nabla\partial_i p(x), w\rangle = \frac{1}{n}\sum_{j=1}^d \omega_d^{-j} \langle \nabla p(x+\omega_d^j e_i), w\rangle
\end{align*}
By Lemma~\ref{lem:fast_grad}, for each $j\in [d]$, $\langle \nabla p(x+\omega_d^j e_i), w\rangle$ can be computed in $O(d\mathcal{T}_O)$ time. Therefore, $\langle \nabla\partial_i p(x), w\rangle$ can be computed in $O(d^2\mathcal{T}_O)$ time.

Consequently, for each coordinate $i\in [w]$, $\left(\nabla^2 (-\ln p(x))w\right)_i$ can be computed in $O(d^2\mathcal{T}_O)$ time. And thus, $\nabla^2 (-\ln p(x))w$ can be computed in $O(nd^2\mathcal{T}_O)$ time.
\end{proof}
\section{List of hyperbolic polynomials} 
\label{sec:list_of_hp}

In this section, we present a list of hyperbolic polynomials related to LP and SDP.

\begin{table*}[!ht]
\centering
\begin{tabular}{|l|l|l|l|l|l|}
\hline
 $p(x)$ & degree &$\nabla(-\ln p(x))[w]$ & $\nabla^2(-\ln p(x))[w]$ & Comment\\ \hline
    $\det(X)$ & $\sqrt{n}$ & Fact~\ref{fac:p1_gradient} & Fact~\ref{fac:p1_hessian} & SDP \\ \hline
    $\prod_{i=1}^n x_i$ & $n$ & Fact~\ref{fac:p2_gradient} & Fact~\ref{fac:p2_hessian} & LP  \\ \hline
    $x_n^2-\sum_{i=1}^{n-1}x_i^2$ & 2 & Fact~\ref{fac:p3_gradient} &  Fact~\ref{fac:p3_hessian} & \\ \hline
    $\det(\sum_{i=1}^n x_i A_i), A_i \in \R^{d\times d}$ & $d$ & Fact~\ref{fac:p4_gradient} & Fact~\ref{fac:p4_hessian} & LP+SDP\\ \hline
\end{tabular}
\caption{Summary of the degree/gradient/Hessian of several hyperbolic polynomials. Note that $\nabla(-\ln p(x))[w] = \langle \nabla(-\ln p(x)), w\rangle$}
\end{table*}

\begin{fact}[Gradient]\label{fac:p1_gradient}
Let $p(X)=\det(X)$, then $\nabla(-\ln p(x))[W]$ is
\begin{align*}
    -\tr(X^{-1}W)
\end{align*}
\end{fact}

\begin{fact}[Hessian]\label{fac:p1_hessian}
\begin{align*}
    X^{-1}WX^{-1}
\end{align*}
\end{fact}

\begin{fact}[Gradient]\label{fac:p2_gradient}
\begin{align*}
    -\sum_{i=1}^n \frac{w_i}{x_i}
\end{align*}

\end{fact}

\begin{fact}[Hessian]\label{fac:p2_hessian}
\begin{align*}
    \begin{bmatrix}\frac{w_1}{x_1^2}&\frac{w_2}{x_2^2}&\cdots& \frac{w_n}{x_n^2}\end{bmatrix}^\top 
\end{align*}
\end{fact}

\begin{fact}[Gradient]\label{fac:p3_gradient}
Let $p(x) = \prod_{i=1}^n x_i$, then $\langle \nabla ( - \ln p(x) ) , w \rangle$ is
\begin{align*}
    2p(x)^{-1}\cdot (x_1w_1+\cdots + x_{n-1}w_{n-1}-x_nw_n)
\end{align*}
\end{fact}

\begin{fact}[Hessian]\label{fac:p3_hessian}
Let $p(x) = \prod_{i=1}^n x_i$, then $\nabla^2 ( - \ln p(x) ) [w]$ is
\begin{align*}
    (\nabla^2(-\ln p(x) & )[w])_i  = \\
    &
    \begin{cases}
        \frac{1}{ p(x)^2 }
         \cdot ( +4x_i\cdot \Gamma_n+2p(x)w_i ),  i\in [n-1];\\
        \frac{1}{p(x)^2} \cdot ( -4x_n \cdot \Gamma_n + 2p(x)w_n ), i=n,
    \end{cases}
\end{align*}
where $\Gamma_n := \sum_{j=1}^{n-1} x_jw_j-x_nw_n $. 
\end{fact}

\begin{fact}[Gradient]\label{fac:p4_gradient}
Let $p(x) = \det( \sum_{i=1}^n x_i A_i )$, then $\langle \nabla ( - \ln p(x) ) , w \rangle$ is
\begin{align*}
    -\tr[ (x_1A_1+\cdots x_nA_n)^{-1}(w_1A_1+\cdots + w_nA_n)]
\end{align*}
\end{fact}

\begin{fact}[Hessian]\label{fac:p4_hessian}
Let $p(x) = \det( \sum_{i=1}^n x_i A_i ), X = \sum_{i=1}^n x_i A_i, W = \sum_{i=1}^n w_i A_i$, then $\nabla^2 ( - \ln p(x) )[w]$ is
\begin{align*}
    \nabla^2   p(x)[w]=
    \begin{bmatrix}
    \tr[X_1] & \tr[X_2] & \cdots & \tr[X_n]
    \end{bmatrix}^\top,
\end{align*}
where $X_i := X^{-1}A_i X^{-1} W$.
\end{fact}
\begin{proof}
Let $g(X)=\log \det(X)$, for $X=\sum_{i=1}^n x_i A_i$. We know that the second derivative of $g(X)$ is given by
\begin{align*}
    \nabla^2 g(X)[Y,Z]=\tr[X^{-1}YX^{-1}Z],
\end{align*}
where $Y,Z\in \R^{n\times n}$. Thus, for $i,j\in [n]\times [n]$,
\begin{align*}
    \nabla^2 p(x)_{i,j}=\tr[X^{-1}A_iX^{-1}A_j].
\end{align*}
Therefore, 
\begin{align*}
    \nabla^2 p(x)[w]= \begin{bmatrix}
    \tr[X_1] & \tr[X_2] & \cdots & \tr[X_n]
    \end{bmatrix},
\end{align*}
where $X_i := X^{-1}A_i X^{-1} (\sum_{i=1}^n w_i A_i)$. 
\end{proof}

\begin{fact}
Let $p(x) = x_1^2-\sum_{i=2}^{n}x_i^2$. The hyperbolic cone of $p(x)$ is
\begin{align*}
    \Lambda_+ = \left\{x\in \R^n: x_1 \geq \sqrt{\sum_{i=2}^n x_i^2}\right\}.
\end{align*}
Furthermore, $\Lambda_+^* = \Lambda_+$.
\end{fact}

\begin{proof}
We first have $e=\begin{bmatrix} 1 & 0 & \cdots & 0 \end{bmatrix}^\top$. Then, we have
\begin{align*}
    p(x-te)=(x_1-t)^2 - \sum_{i=2}^n x_i^2=0,
\end{align*}
which gives
\begin{align*}
    t_{1,2} = x_1\pm \sqrt{\sum_{i=2}^n x_i^2}. 
\end{align*}
Hence, for $x\in \Lambda_+$, we have $t_{1,2}\geq 0$. It follows that $x_1 \geq \sqrt{\sum_{i=2}^n x_i^2}$. Therefore,
\begin{align*}
    \Lambda_+ = \left\{x\in \R^n: x_1 \geq \sqrt{\sum_{i=2}^n x_i^2}\right\}.
\end{align*}

For $s\in \Lambda_+^*$, we have 
\begin{align*}
    \langle s, x\rangle \geq 0 \quad \forall~x\in \Lambda_+,
\end{align*}
That is,
\begin{align*}
    s_1 \geq -\sum_{i=2}^n s_i \frac{x_i}{x_1}.
\end{align*}
Since $x\in \Lambda_+$, 
\begin{align*}
    \sum_{i=2}^n \frac{x_i^2}{x_1^2}\leq 1,
\end{align*}
we can equivalently write the condition of $s\in \Lambda_+^*$ as
\begin{align*}
    s_1 \geq -\sum_{i=2}^n s_i v_i \quad \forall ~\|v\|_2\leq 1. 
\end{align*}
We know that
\begin{align*}
    \RHS \leq \sqrt{\sum_{i=2}^n s_i^2},
\end{align*}
and the equality holds if $v_i = -s_i/\sqrt{\sum_{i=2}^n s_i^2}$ for $2\leq i\leq n$. Therefore, $s\in \Lambda_+^*$ if and only if
\begin{align*}
    s_1 \geq \sqrt{\sum_{i=2}^n s_i^2}.
\end{align*}
That is,
\begin{align*}
    \Lambda_+^* = \Lambda_+.
\end{align*}
\end{proof}

\begin{fact}
Let $p(x)=\det(\sum_{i=1}^n x_i A_i)$ for $A_i \in \R^{d\times d}$ are real symmetric matrices and there exists an $e\in \R^n$ such that $\sum_{i=1}^n e_i A_i = I_d$. The hyperbolic cone of $p(x)$ is spectrahedral (Definition~\ref{def:spec_hed_cone}), i.e.,
\begin{align*}
    \Lambda_+ = \left\{x\in \R^n: \sum_{i=1}^n x_i A_i  \succeq 0\right\}.
\end{align*}
\end{fact}
\begin{proof}
Let $E:= \sum_{i=1}^n e_i A_i$. By definition, $E\succ 0$. Consider $p(x-te)=\det(\sum_{i=1}^n x_i A_i - tE)=$
\begin{align*}
    p(x-te)=&~ \det\left(\sum_{i=1}^n x_i A_i - tE\right)\\
    = &~ \det\left(E^{-1/2}\left(\sum_{i=1}^n x_i A_i\right) E^{-1/2} - tI\right)\det(E)
\end{align*}
Since $A_i$ and $E$ are symmetric, $p(x-te)$ has only real roots. And the roots, or the eigenvalues of $E^{-1/2}\left(\sum_{i=1}^n x_i A_i\right) E^{-1/2}$, greater or equal to 0 if and only if $\sum_{i=1}^n x_i A_i$ is a PSD matrix.
\end{proof}

\begin{fact}
The dual of spectrahedral cone may not be spectrahedral.
\end{fact}

\ifdefined\isarxiv
\bibliographystyle{alpha}
\bibliography{ref}

\newcommand{\etalchar}[1]{$^{#1}$}
\begin{thebibliography}{SWYZ23}

\bibitem[AW21]{aw21}
Josh Alman and Virginia~Vassilevska Williams.
\newblock A refined laser method and faster matrix multiplication.
\newblock In {\em Proceedings of the 2021 ACM-SIAM Symposium on Discrete
  Algorithms (SODA)}, pages 522--539. SIAM, 2021.

\bibitem[Br{\"a}14]{bra14}
Petter Br{\"a}nd{\'e}n.
\newblock Hyperbolicity cones of elementary symmetric polynomials are
  spectrahedral.
\newblock {\em Optimization Letters}, 8(5):1773--1782, 2014.

\bibitem[Bra20]{b20}
Jan van~den Brand.
\newblock A deterministic linear program solver in current matrix
  multiplication time.
\newblock In {\em Proceedings of the Fourteenth Annual ACM-SIAM Symposium on
  Discrete Algorithms (SODA)}, pages 259--278. SIAM, 2020.

\bibitem[CLS19]{cls19}
Michael~B Cohen, Yin~Tat Lee, and Zhao Song.
\newblock Solving linear programs in the current matrix multiplication time.
\newblock In {\em Proceedings of the 51st annual ACM SIGACT symposium on theory
  of computing}, pages 938--942, 2019.

\bibitem[Dan47]{d47}
George~B Dantzig.
\newblock Maximization of a linear function of variables subject to linear
  inequalities.
\newblock {\em Activity analysis of production and allocation}, 13:339--347,
  1947.

\bibitem[DLY21]{dly21}
Sally Dong, Yin~Tat Lee, and Guanghao Ye.
\newblock A nearly-linear time algorithm for linear programs with small
  treewidth: A multiscale representation of robust central path.
\newblock In {\em Proceedings of the 53rd Annual ACM SIGACT Symposium on Theory
  of Computing (STOC)}, 2021.

\bibitem[G{\aa}r51]{g51}
Lars G{\aa}rding.
\newblock Linear hyperbolic partial differential equations with constant
  coefficients.
\newblock {\em Acta Mathematica}, 85:1--62, 1951.

\bibitem[G{\aa}r59]{g59}
Lars G{\aa}rding.
\newblock An inequality for hyperbolic polynomials.
\newblock {\em Journal of Mathematics and Mechanics}, pages 957--965, 1959.

\bibitem[GS22]{gs22}
Yuzhou Gu and Zhao Song.
\newblock A faster small treewidth sdp solver.
\newblock {\em arXiv preprint arXiv:2211.06033}, 2022.

\bibitem[G{\"u}l97]{gul97}
Osman G{\"u}ler.
\newblock Hyperbolic polynomials and interior point methods for convex
  programming.
\newblock {\em Mathematics of Operations Research}, 22(2):350--377, 1997.

\bibitem[HJS{\etalchar{+}}22]{hjstz22}
Baihe Huang, Shunhua Jiang, Zhao Song, Runzhou Tao, and Ruizhe Zhang.
\newblock Solving sdp faster: A robust ipm framework and efficient
  implementation.
\newblock In {\em FOCS}, 2022.

\bibitem[HV07]{hv07}
J~William Helton and Victor Vinnikov.
\newblock Linear matrix inequality representation of sets.
\newblock {\em Communications on Pure and Applied Mathematics: A Journal Issued
  by the Courant Institute of Mathematical Sciences}, 60(5):654--674, 2007.

\bibitem[Jia21]{j21}
Haotian Jiang.
\newblock Minimizing convex functions with integral minimizers.
\newblock In {\em Proceedings of the 2021 ACM-SIAM Symposium on Discrete
  Algorithms}, 2021.

\bibitem[Jia22]{j22}
Haotian Jiang.
\newblock Minimizing convex functions with rational minimizers.
\newblock {\em ACM Journal of the ACM (JACM)}, 2022.

\bibitem[JKL{\etalchar{+}}20]{jklps20}
Haotian Jiang, Tarun Kathuria, Yin~Tat Lee, Swati Padmanabhan, and Zhao Song.
\newblock A faster interior point method for semidefinite programming.
\newblock 2020.

\bibitem[JLSW20]{jlsw20}
Haotian Jiang, Yin~Tat Lee, Zhao Song, and Sam Chiu-wai Wong.
\newblock An improved cutting plane method for convex optimization,
  convex-concave games, and its applications.
\newblock In {\em Proceedings of the 52nd Annual ACM SIGACT Symposium on Theory
  of Computing}, pages 944--953, 2020.

\bibitem[JLSZ23]{jlsz23}
Haotian Jiang, Yin~Tat Lee, Zhao Song, and Lichen Zhang.
\newblock Convex minimization with integral minima in $\widetilde {O}(n^4)$
  time.
\newblock {\em arXiv preprint arXiv:2304.03426}, 2023.

\bibitem[JSWZ21]{jswz21}
Shunhua Jiang, Zhao Song, Omri Weinstein, and Hengjie Zhang.
\newblock Faster dynamic matrix inverse for faster lps.
\newblock In {\em Proceedings of the 53rd Annual ACM SIGACT Symposium on Theory
  of Computing (STOC)}, 2021.

\bibitem[Kar84]{k84}
N.~Karmarkar.
\newblock A new polynomial-time algorithm for linear programming.
\newblock In {\em Proceedings of the Sixteenth Annual ACM Symposium on Theory
  of Computing}, STOC '84, 1984.

\bibitem[Kha80]{k80}
Leonid~G Khachiyan.
\newblock Polynomial algorithms in linear programming.
\newblock {\em USSR Computational Mathematics and Mathematical Physics},
  20(1):53--72, 1980.

\bibitem[Lax57]{lax57}
Peter~D Lax.
\newblock Differential equations, difference equations and matrix theory.
\newblock Technical report, New York Univ., New York. Atomic Energy Commission
  Computing and Applied, 1957.

\bibitem[LPR05]{lpr05}
Adrian Lewis, Pablo Parrilo, and Motakuri Ramana.
\newblock The lax conjecture is true.
\newblock {\em Proceedings of the American Mathematical Society},
  133(9):2495--2499, 2005.

\bibitem[LS14]{ls14}
Yin~Tat Lee and Aaron Sidford.
\newblock Path finding methods for linear programming: Solving linear programs
  in ${O}(\sqrt{rank})$ iterations and faster algorithms for maximum flow.
\newblock In {\em 2014 IEEE 55th Annual Symposium on Foundations of Computer
  Science}, pages 424--433. IEEE, 2014.

\bibitem[LS19]{ls19}
Yin~Tat Lee and Aaron Sidford.
\newblock Solving linear programs with sqrt (rank) linear system solves.
\newblock {\em arXiv preprint arXiv:1910.08033}, 2019.

\bibitem[LSW15]{lsw15}
Yin~Tat Lee, Aaron Sidford, and Sam Chiu-wai Wong.
\newblock A faster cutting plane method and its implications for combinatorial
  and convex optimization.
\newblock In {\em 2015 IEEE 56th Annual Symposium on Foundations of Computer
  Science}, pages 1049--1065. IEEE, 2015.

\bibitem[LSZ19]{lsz19}
Yin~Tat Lee, Zhao Song, and Qiuyi Zhang.
\newblock Solving empirical risk minimization in the current matrix
  multiplication time.
\newblock In {\em COLT}, 2019.

\bibitem[LSZ{\etalchar{+}}20]{lsz+20}
S~Cliff Liu, Zhao Song, Hengjie Zhang, Lichen Zhang, and Tianyi Zhou.
\newblock Space-efficient interior point method, with applications to linear
  programming and maximum weight bipartite matching.
\newblock {\em arXiv e-prints}, pages arXiv--2009, 2020.

\bibitem[Pem12]{pem12}
Robin Pemantle.
\newblock Hyperbolicity and stable polynomials in combinatorics and
  probability.
\newblock {\em arXiv preprint arXiv:1210.3231}, 2012.

\bibitem[QSZZ23]{qszz23}
Lianke Qin, Zhao Song, Lichen Zhang, and Danyang Zhuo.
\newblock An online and unified algorithm for projection matrix vector
  multiplication with application to empirical risk minimization.
\newblock In {\em AISTATS}, 2023.

\bibitem[Ren88]{r88}
James Renegar.
\newblock A polynomial-time algorithm, based on {N}ewton's method, for linear
  programming.
\newblock {\em Math. Programming}, 1988.

\bibitem[Ren95]{r95}
James Renegar.
\newblock Linear programming, complexity theory and elementary functional
  analysis.
\newblock {\em Math. Programming}, 1995.

\bibitem[Ren04]{ren04}
James Renegar.
\newblock Hyperbolic programs, and their derivative relaxations.
\newblock Technical report, Cornell University Operations Research and
  Industrial Engineering, 2004.

\bibitem[RS14]{rs14}
James Renegar and Mutiara Sondjaja.
\newblock A polynomial-time affine-scaling method for semidefinite and
  hyperbolic programming.
\newblock {\em arXiv preprint arXiv:1410.6734}, 2014.

\bibitem[Sau19]{s19}
James Saunderson.
\newblock Certifying polynomial nonnegativity via hyperbolic optimization.
\newblock {\em SIAM Journal on Applied Algebra and Geometry}, 3(4):661--690,
  2019.

\bibitem[SWYZ21]{swyz21}
Zhao Song, David Woodruff, Zheng Yu, and Lichen Zhang.
\newblock Fast sketching of polynomial kernels of polynomial degree.
\newblock In {\em International Conference on Machine Learning}, pages
  9812--9823. PMLR, 2021.

\bibitem[SWYZ23]{swyz23}
Zhao Song, Yitan Wang, Zheng Yu, and Lichen Zhang.
\newblock Sketching for first order method: efficient algorithm for
  low-bandwidth channel and vulnerability.
\newblock In {\em International Conference on Machine Learning (ICML)}, 2023.

\bibitem[SY21]{sy21}
Zhao Song and Zheng Yu.
\newblock Oblivious sketching-based central path method for solving linear
  programming problems.
\newblock In {\em 38th International Conference on Machine Learning (ICML)},
  2021.

\bibitem[SYYZ23]{syyz23}
Zhao Song, Xin Yang, Yuanyuan Yang, and Lichen Zhang.
\newblock Sketching meets differential privacy: fast algorithm for dynamic
  kronecker projection maintenance.
\newblock In {\em International Conference on Machine Learning (ICML)}, 2023.

\bibitem[SYZ21]{syz21}
Zhao Song, Shuo Yang, and Ruizhe Zhang.
\newblock Does preprocessing help training over-parameterized neural networks?
\newblock {\em Advances in Neural Information Processing Systems},
  34:22890--22904, 2021.

\bibitem[SZZ21]{szz21}
Zhao Song, Lichen Zhang, and Ruizhe Zhang.
\newblock Training multi-layer over-parametrized neural network in subquadratic
  time.
\newblock {\em arXiv preprint arXiv:2112.07628}, 2021.

\bibitem[Vai89a]{v89}
Pravin~M Vaidya.
\newblock A new algorithm for minimizing convex functions over convex sets.
\newblock In {\em 30th Annual IEEE Symposium on Foundations of Computer Science
  (FOCS)}, 1989.

\bibitem[Vai89b]{v89_lp}
Pravin~M Vaidya.
\newblock Speeding-up linear programming using fast matrix multiplication.
\newblock In {\em 30th Annual Symposium on Foundations of Computer Science
  (FOCS)}, 1989.

\bibitem[Vin12]{vin12}
Victor Vinnikov.
\newblock Lmi representations of convex semialgebraic sets and determinantal
  representations of algebraic hypersurfaces: past, present, and future.
\newblock In {\em Mathematical methods in systems, optimization, and control},
  pages 325--349. Springer, 2012.

\bibitem[Wag11]{wag11}
David Wagner.
\newblock Multivariate stable polynomials: theory and applications.
\newblock {\em Bulletin of the American Mathematical Society}, 48(1):53--84,
  2011.

\bibitem[Ye21]{y21}
Guanghao Ye.
\newblock Fast algorithm for solving structured convex programs.
\newblock Master's thesis, The University of Washington, 2021.

\end{thebibliography}

\else

\fi

\end{document}